\documentclass[DIV=15,11pt]{scrartcl}

\usepackage[utf8]{inputenc}
\usepackage[english]{babel}
\DeclareUnicodeCharacter{0301}{\'{e}}
\usepackage{amsmath,amssymb,amsfonts,amsthm}
\usepackage{stmaryrd}
\usepackage{graphicx}
\usepackage{wrapfig}
\usepackage{subcaption}
\usepackage{multirow}
\usepackage{booktabs}
\usepackage{tabularx}
\usepackage{enumerate}

\usepackage{algorithm}
\usepackage{algpseudocode}

\usepackage{tikz}
\usetikzlibrary{calc,math,quotes,angles,positioning}

\usepackage[colorlinks=true,urlcolor=blue,citecolor=blue,linkcolor=blue,bookmarks=true]{hyperref}

\usepackage{amsthm}
\theoremstyle{definition}

\usepackage{tcolorbox}
\usepackage{underscore}
\usepackage[autostyle]{csquotes}
\newtheorem{lemma}{Lemma}[section]

\newtheorem{theorem}{Theorem}[section]

\graphicspath{images/}

\DeclareMathOperator{\supp}{supp} 

\usepackage{algpseudocode}

\newcommand{\vertiii}[1]{{\left\vert\kern-0.25ex\left\vert\kern-0.25ex\left\vert #1 
    \right\vert\kern-0.25ex\right\vert\kern-0.25ex\right\vert}}
    
\begin{document}

\title{A non-negativity-preserving cut-cell discontinuous Galerkin method for the diffusive wave equation}

\author{P. Manorost, P. Bastian\\
Interdisciplinary Center for Scientific Computing\\
Heidelberg University\\
Im Neuenheimer Feld 205, D-69120 Heidelberg\\
\texttt{peter.bastian@iwr.uni-heidelberg.de}}

\maketitle
\begin{abstract}
A non-negativity-preserving cut-cell discontinuous Galerkin method for the degenerate parabolic diffusive wave approximation of the shallow water equation is presented. The method can handle continuous and discontinuous bathymmetry as well as general triangular meshes. It is complemented by a finite volume method on Delauney triangulations which is also shown to be non-negativity preserving. Both methods feature an upwind flux and can handle Manning's and Chezy's friction law. By numerical experiment we demonstrate the discontinuous Galerkin method to be fully second-order accurate for the Barenblatt analytical solution on an inclined plane. In constrast, the finite volume method is only first-order accurate. Further numerical experiments show that three to four mesh refinements are needed for the finite volume method to match the solution of the discontinuous Galerkin method. 
\end{abstract}

\section{Introduction}

To understand and predict surface flows is of enormous importance to society for minimizing risks and impact of flooding events, e.g. in coastal areas \cite{DIETRICH201145}, urban areas \cite{MIGNOT2006186} and flood plains  \cite{bauer2006regional} or, coupled with subsurface flow, to address questions of water availability and quality \cite{Wood2011,PaniconiPutti2015}. 
Only recently, e.g., hydrology in high mountain catchments has been identified as a region playing a pivotal role in the hydrological cycle \cite{VINCENT201959,Thornton2022,vantiel2024,SOHEB2024131063}, necessitating accurate predictions.

Physics-based modelling of surface flows is often based on the shallow water equations (SWE), also termed dynamic wave equation, which is a depth-integrated simplification of the incompressible Navier-Stokes equations derived under the assumption of negligable vertical velocity and horizontal velocity components not deviating much from their depth-average \cite{tan1992shallow}. Neglegting internal friction, the SWE comprise a nonlinear first-order hyperbolic system of partial differential equations. In this paper we are concerned with a further simplification of the SWE termed diffusive wave equation (DWE) or zero inertia equation (ZI), which is obtained from the SWE under the assumption that gravity and bed friction terms are dominating \cite{ponce1977shallow,BOLSTER2002221,HUNTER2007208,savant2019urban}. The DWE is a scalar nonlinear degenerate parabolic partial differential equation which is a good approximation of the SWE in case of slowly evolving flows, for comparisons and validity we refer to \cite{HUNTER2007208,CEA201088,Martins2017,COSTABILE2017141,Arico2018,CAVIEDESVOULLIEME2020124663}. 

As a scalar equation the DWE is deemed to be numerically less costly to simulate than the SWE. However, several studies have found its numerical solution to be more costly \cite{CEA201088,CAVIEDESVOULLIEME2020124663}. This was a consequence of using explicit time stepping schemes, a fact well-known in numerical analysis for parabolic partial differential equations. A formal derivation of the stability constraint is given in \cite{hunter2005adaptive}. Using implicit time-stepping the time step restriction is avoided at the cost of solving large scale nonlinear algebraic systems \cite{SZYMKIEWICZ2012165,fernandez20162d}. When using efficient solvers relatively large times steps can be taken \cite{SOHEB2024131063}. 

As spatial discretization the Finite Volume method (FVM) is hugely popular \cite{hunter2005adaptive,CEA201088,wang2011positivity,caviedes2018cellular,CAVIEDESVOULLIEME2020124663}, in particular on rectangular grids. This is in part due to the abundant availability of digital elevation data in raster form. 
Finite element discretizations on triangular unstructured meshes have also been presented by several authors.
\cite{szymkiewicz2012simulation} presents a modified lowest-order conforming finite element method where the modification refers to the mass matrix.
\cite{savant2019urban} presents a lowest-order conforming finite element method stabilized with a Galerkin least squares formulation and featuring adaptive mesh refinement.
\cite{de2019discontinuous} develops a control volume finite element method for the DWE and
\cite{santillana2010numerical} uses a conforming finite element method with mass lumping in one spatial dimension. A local discontinuous Galerkin finite method for the DWE has been presented in \cite{santillana2010local} for continuous bathymmetry and demonstrating a wetting front but no drying front. Cut-cell (or unfitted) finite elements, including discontinuous Galerkin methods, have been developed for solving partial differential equations in complex domains \cite{HANSBO20025537,Bastian2009,Burman2015,BADIA2018533}. In our application the complex domain is defined dynamically by the support of the water height which has the wet/dry front as boundary. Cut-cell methods for free boundary problems have been presented in \cite{schneiders2013accurate,Heimann2013,liu2021parallel}. To our knowledge the application to the DWE with wet/dry front is new.

In this work we combine the discontinous Galerkin and cut-cell approach for the efficient and accurate numerical solution of the DWE. The main contributions of the paper are:
\begin{itemize}
\item Presentation of a fully second-order accurate, locally conservative discontinuous Galerkin (DG) finite element method (FEM) for the DWE,
\item that is able to handle continuous and discontinuous representations of bathymmetry on unstructured triangular meshes. Moreover, the method allows handling of the norm of the gradient of free surface elevation in the friction law.
\item Implicit time-stepping using Newton's method to solve the nonlinear algebraic systems avoids a time step restriction for stability reasons.
\item A new cut-cell approach ensures non-negative water height in the presence of wet/dry fronts and second-order accuracy up to the free boundary is demonstrated numerically.
\item Numerical examples compare the DG method with a Voronoi finite volume method (VFVM) equipped with a four-point flux for the full DWE model. In contrast to the finite volume method the discontinuous Galerkin method allows for more general meshes including local mesh refinement. 
\end{itemize}

The rest of the paper is organized as follows. In the next section, the mathematical model is introduced together with its basic properties. In section three we present the Voronoi finite volume scheme and the cut-cell discontinuous Galerkin finite element method. Numerical results are presented in section four, demonstrating the behaviour of the methods on three test problems of increasing difficulty. The results of the paper are summarized in section five.

\section{Mathematical model} \label{sec:math_model}

Let $\Omega$ be a polygonal domain in $\mathbb{R}^2$ and $\Sigma=(t_0,t_0+T)$ a time interval. The diffusive wave equation (DWE) reads \cite{ALONSOSANTILLANADAWSON2008,cea2010experimental,costabile2017performances,savant2019urban,bogelein2021existence}: 
\begin{subequations} \label{model_problem_peter}
\begin{align}
\partial_t u(x,t) + \nabla \cdot q(x,t) &= f(x,t), &&(x,t) \in \Omega \times \Sigma \label{eq:conservation}\\
q(x,t) &= - K(x) \frac{H(x,t)^\alpha}{G(x,t)^{1-\gamma}} \ \nabla u(x,t) \label{eq:vol_flux}  \\
u(x,t) &= g(x,t), &&(x,t) \in \Gamma_D\times\Sigma \label{eq:DWE-Dirichlet} \\
q(x,t) \cdot n(x) &= j(x,t), &&(x,t) \in \Gamma_N\times\Sigma \label{eq:DWE-Neumann}\\
u(x,t_0) &= u_0(x) && x\in\Omega, \label{eq:DWE-Initial}
\end{align}
\end{subequations}
where 
$u(x,t)$ is the free water surface elevation $[m]$,
$q(x,t)$ is the volumetric flux $[m^2/s]$,
\begin{equation*}
H(x,t) = u(x,t)-b(x) \qquad\text{and}\qquad G(x,t) = \|\nabla u(x,t)\|
\end{equation*}
are the height of the water column over the given time-independent bathymmetry (or land elevation) $b(x)$ and 
the Euclidean norm of the free surface elevation gradient.
$K(x)$ is a scalar function which, together with the parameters 
$\alpha\in (0,2]$ and $\gamma\in (0,1]$ models different friction laws.
Most common are Manning's friction law given by 
\begin{equation*}
K(x) = 1/n, \qquad \alpha = 5/3, \qquad \gamma = 1/2, 
\end{equation*}
with $n$ being Manning's number and Chezy's formula, where 
\begin{equation*}
K(x) = C, \qquad \alpha = 3/2, \qquad \gamma = 1/2. 
\end{equation*}
Setting $\alpha=1$, $\gamma=1$ and $K(x)$ the hydraulic conductivity models depth-averaged groundwater flow in an unconfined aquifer.
Volumetric flux $q$ can also be written in terms of $q(x,t) = H(x,t) V(x,t)$ with the
flow velocity $V(x,t)=-K(x)H(x,t)^{\alpha-1}G(x,t)^{\gamma-1}\nabla u(x,t)$ in $[m/s]$. 
Equations \eqref{eq:DWE-Dirichlet}, \eqref{eq:DWE-Neumann} are the Dirichlet and Neumann boundary conditions (with $n(x)$ the unit outer normal vector at the boundary) and \eqref{eq:DWE-Initial} is the initial condition.

Equation \eqref{model_problem_peter} is of doubly degenerate parabolic type. Existence, regularity, uniqueness and non-negativity for a constant bathymmetry are studied in \cite{ALONSOSANTILLANADAWSON2008}. A numerical approach to study properties of the DWE is presented in \cite{santillana2010numerical}. Existence of non-negative weak solutions for non-constant bathymmetry is shown in \cite{bogelein2021existence} and local boundedness is shown in \cite{singer2019local}. 
Non-negativity of the water height $H(x,t) = u(x,t) - b(x)$ is a very important structural property of the solution of the DWE, which needs to be maintained by any numerical solution. 


\section{Numerical schemes} \label{math section}

In this section we describe two numerical schemes: A finite volume method based on Voronoi cells and a cut-cell discontinuous Galerkin method.
Both schemes are based on a partitioning of the domain $\Omega\subset\mathbb{R}^2$ into a set of open, nonoverlapping cells $\mathcal{C}_h = \{ C_1, \ldots, C_M\}$ to be specified below. Each cell $C\in\mathcal{C}_h$ is itself polygonal with a center $x_C\in C$. The cells $\mathcal{C}_h^D = \{ C\in\mathcal{C}_h : \partial C \cap \Gamma_D \neq \emptyset\}$ are those on the Dirichlet boundary and by $\mathcal{C}_h^{iN} = \mathcal{C}_h\setminus\mathcal{C}_h^D$ we denote the interior and Neumann boundary cells. 

We call a line segment $F = \partial C \cap \Gamma_D$ a Dirichlet boundary edge and $F = \partial C \cap \Gamma_N$ a Neumann boundary edge. For any given boundary edge $F$, we denote by $C_F^-$ its associated cell, by $x_F\in\partial\Omega$ its center and by $n_F$ its associated unit normal vector which coincides with the unit outer normal vector at $x_F$. Boundary edges may be subdivided into individual straight line segments to ensure a unique normal direction for each segment.

Given two cells $C,C'\in\mathcal{C}_h$ their intersection $F = \partial C \cap \partial C'$ is called an interior edge. For any interior edge $F$ we denote by $C^-_F, C^+_F$ the two associated cells with respect to a chosen unit normal vector $n_F$ to the edge $F$ pointing from $C^-_F$ to $C^+_F$. The center of $F$ is denoted by $x_F$. The set of all interior edges is denoted by $\mathcal{F}_h^i$, the set of all Neumann boundary 
edges is $\mathcal{F}_h^N$, the set of all Dirichlet boundary edges is $\mathcal{F}_h^D$
and the set of all interior and Dirichlet boundary edges is $\mathcal{F}_h^{iD}$.

The temporal domain $\Sigma$ is discretized into intervalls by $t_0 = t^0 < t^1 < \ldots < t^N = t_0+T\}$, where $N$ is a positive integer. The time step size is given by $\Delta t^n = t^{n} - t^{n-1}$.

\subsection{Voronoi finite volume method}\label{sec:VFVM}

\begin{figure}
\begin{center}
\begin{tikzpicture}[scale=0.9]
\tikzstyle{every node}=[font=\scriptsize]

\coordinate (v1) at (0,1);
\coordinate (v2) at (2,0);
\coordinate (v3) at (2,2);
\coordinate (v4) at (3.8,1.8);
\coordinate (v5) at (0,3);
\coordinate (v6) at (3.5,3.5);
\coordinate (v7) at (2,5);
\coordinate (v8) at (4.5,0);
\coordinate (v9) at (6,2.5);
\coordinate (v10) at (5,5);

\draw[lightgray,thin] (v3) -- (v1);
\draw[lightgray,thin] (v3) -- (v2);
\draw[lightgray,thin] (v3) -- (v4);
\draw[lightgray,thin] (v3) -- (v5);
\draw[lightgray,thin] (v3) -- (v6);
\draw[lightgray,thin] (v3) -- (v7);
\draw[lightgray,thin] (v1) -- (v2);

\draw[lightgray,thin] (v2) -- (v4);
\draw[lightgray,thin] (v4) -- (v6);
\draw[lightgray,thin] (v6) -- (v7);
\draw[lightgray,thin] (v7) -- (v5);
\draw[lightgray,thin] (v5) -- (v1);
\draw[lightgray,thin] (v2) -- (v8);
\draw[lightgray,thin] (v4) -- (v9);
\draw[lightgray,thin] (v4) -- (v8);
\draw[lightgray,thin] (v8) -- (v9);
\draw[lightgray,thin] (v6) -- (v9);
\draw[lightgray,thin] (v6) -- (v10);
\draw[lightgray,thin] (v10) -- (v9);
\draw[lightgray,thin] (v7) -- (v10);

\coordinate (m23) at ($(v2)!0.5!(v3)$);
\coordinate (t1) at ($(m23)!0.71!90:(v3)$);
\coordinate (t2) at ($(m23)!0.9!270:(v3)$);
\coordinate (m24) at ($(v2)!0.5!(v4)$);
\coordinate (t3) at ($(m24)!0.4!270:(v4)$);
\coordinate (m13) at ($(v1)!0.5!(v3)$);
\coordinate (t4) at ($(m13)!0.55!90:(v3)$);
\coordinate (m53) at ($(v5)!0.5!(v3)$);
\coordinate (t5) at ($(m53)!1.0!90:(v3)$);
\coordinate (m73) at ($(v7)!0.5!(v3)$);
\coordinate (t6) at ($(m73)!0.02!90:(v3)$);
\coordinate (m63) at ($(v6)!0.5!(v3)$);
\coordinate (t7) at ($(m63)!0.5!90:(v3)$);
\coordinate (m12) at ($(v1)!0.5!(v2)$);
\coordinate (m15) at ($(v1)!0.5!(v5)$);
\coordinate (m57) at ($(v5)!0.5!(v7)$);
\coordinate (m28) at ($(v2)!0.5!(v8)$);
\coordinate (m48) at ($(v4)!0.5!(v8)$);
\coordinate (t8) at ($(m48)!1.1!90:(v8)$);
\coordinate (m49) at ($(v4)!0.5!(v9)$);
\coordinate (t9) at ($(m49)!0.5!90:(v9)$);
\coordinate (m610) at ($(v6)!0.5!(v10)$);
\coordinate (t10) at ($(m610)!1.05!270:(v10)$);
\coordinate (t11) at ($(m610)!1.0!90:(v10)$);
\coordinate (m89) at ($(v8)!0.5!(v9)$);
\coordinate (m910) at ($(v9)!0.5!(v10)$);
\coordinate (m710) at ($(v7)!0.5!(v10)$);

\draw[ultra thick] (t1) -- (t2);
\draw[ultra thick] (t2) -- (t7);
\draw[ultra thick] (t6) -- (t5);
\draw[ultra thick] (t5) -- (t4);
\draw[ultra thick] (t1) -- (t4);
\draw[ultra thick] (m12) -- (t1);
\draw[ultra thick] (m15) -- (t4);
\draw[ultra thick] (m57) -- (t5);
\draw[ultra thick] (t2) -- (t3);
\draw[ultra thick] (m28) -- (t3);
\draw[ultra thick] (t3) -- (t8);
\draw[ultra thick] (t8) -- (t9);
\draw[ultra thick] (t7) -- (t9);
\draw[ultra thick] (t9) -- (t10);
\draw[ultra thick] (t6) -- (t11);
\draw[ultra thick] (t10) -- (t11);
\draw[ultra thick] (m89) -- (t8);
\draw[ultra thick] (m910) -- (t10);
\draw[ultra thick] (m710) -- (t11);
\draw[ultra thick] (v1) -- (v2);
\draw[ultra thick] (v2) -- (v8);
\draw[ultra thick] (v8) -- (v9);
\draw[ultra thick] (v9) -- (v10);
\draw[ultra thick] (v10) -- (v7);
\draw[ultra thick] (v7) -- (v5);
\draw[ultra thick] (v5) -- (v1);

\draw[red,ultra thick] (t7) -- (t6);
\draw[red,thick,->] (m63) -- ($(m63)!0.7!(v6)$) node [sloped,midway,below] {$n_F$} ;
\filldraw [fill=white,draw=red] (m63) circle (3pt); \node[below=3pt] at (m63) {\color{red}$x_F$};
\node[red] at ($(m63)!0.6!330:(t6)$) {$F$};
\filldraw [fill=white,draw=blue] (m63) circle (3pt); \node[below=3pt] at (m63) {\color{red}$x_F$};

\draw[blue,ultra thick] (m89) -- (v9);
\coordinate (xF) at ($(m89)!0.5!(v9)$);
\draw[blue,thick,->] (xF) -- ($(xF)!0.9!270:(v9)$) node [sloped,midway,below] {$n_{F'}$} ;
\filldraw [fill=white,draw=blue] (xF) circle (3pt); \node[above left=-3pt] at (xF) {\color{blue}$x_{F'}$};
\node[blue] at ($(xF)!0.7!320:(v9)$) {$F'$};

\filldraw  [fill=white,draw=black]  (v1) circle (3pt);
\filldraw  [fill=white,draw=black]  (v2) circle (3pt);
\filldraw  [fill=white,draw=black]  (v3) circle (3pt);
\filldraw  [fill=white,draw=black]  (v4) circle (3pt);
\filldraw  [fill=white,draw=black]  (v5) circle (3pt);
\filldraw  [fill=white,draw=black]  (v6) circle (3pt);
\filldraw  [fill=white,draw=black]  (v7) circle (3pt);
\filldraw  [fill=white,draw=black]  (v8) circle (3pt);
\filldraw  [fill=white,draw=black]  (v9) circle (3pt);
\filldraw  [fill=white,draw=black]  (v10) circle (3pt);

\node[font=\small] at (3,-1) {a) Voronoi cells}; 

\node[left] at (v1) {$x_1$};
\node[below] at (v2) {$x_2$};
\node[right] at (v3) {$x_{C^-_F}$};
\node[right] at (v4) {$x_4$};
\node[left] at (v5) {$x_5$};
\node[right] at (v6) {$x_{C^+_F}$};
\node[above] at (v7) {$x_7$};
\node[right] at (v8) {$x_8$};
\node[right] at (v9) {$x_{C_{F'}}$};
\node[right] at (v10) {$x_{10}$};
\end{tikzpicture}
\hfill
\begin{tikzpicture}[scale=0.9]
\tikzstyle{every node}=[font=\scriptsize]

\coordinate (v1) at (0,1);
\coordinate (v2) at (2,0);
\coordinate (v3) at (2,2);
\coordinate (v4) at (3.8,1.8);
\coordinate (v5) at (0,3);
\coordinate (v6) at (3.5,3.5);
\coordinate (v7) at (2,5);
\coordinate (v8) at (4.5,0);
\coordinate (v9) at (6,2.5);
\coordinate (v10) at (5,5);

\draw[ultra thick] (v3) -- (v1);
\draw[ultra thick] (v3) -- (v2);
\draw[ultra thick] (v3) -- (v4);
\draw[ultra thick] (v3) -- (v5);
\draw[ultra thick] (v3) -- (v6);
\draw[ultra thick] (v3) -- (v7);
\draw[ultra thick] (v1) -- (v2);
\draw[ultra thick] (v2) -- (v4);
\draw[ultra thick] (v4) -- (v6);
\draw[ultra thick] (v6) -- (v7);
\draw[ultra thick] (v7) -- (v5);
\draw[ultra thick] (v5) -- (v1);
\draw[ultra thick] (v2) -- (v8);
\draw[ultra thick] (v4) -- (v9);
\draw[ultra thick] (v4) -- (v8);
\draw[ultra thick] (v8) -- (v9);
\draw[ultra thick] (v6) -- (v9);
\draw[ultra thick] (v6) -- (v10);
\draw[ultra thick] (v10) -- (v9);
\draw[ultra thick] (v7) -- (v10);

\coordinate (m23) at ($(v2)!0.5!(v3)$);
\coordinate (t1) at ($(m23)!0.71!90:(v3)$);
\coordinate (t2) at ($(m23)!0.9!270:(v3)$);
\coordinate (m24) at ($(v2)!0.5!(v4)$);
\coordinate (t3) at ($(m24)!0.4!270:(v4)$);
\coordinate (m13) at ($(v1)!0.5!(v3)$);
\coordinate (t4) at ($(m13)!0.55!90:(v3)$);
\coordinate (m53) at ($(v5)!0.5!(v3)$);
\coordinate (t5) at ($(m53)!1.0!90:(v3)$);
\coordinate (m73) at ($(v7)!0.5!(v3)$);
\coordinate (t6) at ($(m73)!0.02!90:(v3)$);
\coordinate (m63) at ($(v6)!0.5!(v3)$);
\coordinate (t7) at ($(m63)!0.5!90:(v3)$);
\coordinate (m12) at ($(v1)!0.5!(v2)$);
\coordinate (m15) at ($(v1)!0.5!(v5)$);
\coordinate (m57) at ($(v5)!0.5!(v7)$);
\coordinate (m28) at ($(v2)!0.5!(v8)$);
\coordinate (m48) at ($(v4)!0.5!(v8)$);
\coordinate (m46) at ($(v4)!0.5!(v6)$);
\coordinate (t8) at ($(m48)!1.1!90:(v8)$);
\coordinate (m49) at ($(v4)!0.5!(v9)$);
\coordinate (t9) at ($(m49)!0.5!90:(v9)$);
\coordinate (m610) at ($(v6)!0.5!(v10)$);
\coordinate (t10) at ($(m610)!1.05!270:(v10)$);
\coordinate (t11) at ($(m610)!1.0!90:(v10)$);
\coordinate (m89) at ($(v8)!0.5!(v9)$);
\coordinate (m910) at ($(v9)!0.5!(v10)$);
\coordinate (m710) at ($(v7)!0.5!(v10)$);


\draw[ultra thick] (v1) -- (v2);
\draw[ultra thick] (v2) -- (v8);
\draw[ultra thick] (v8) -- (v9);
\draw[ultra thick] (v9) -- (v10);
\draw[ultra thick] (v10) -- (v7);
\draw[ultra thick] (v7) -- (v5);
\draw[ultra thick] (v5) -- (v1);

\draw[red,ultra thick] (v4) -- (v6);
\draw[red,thick,->] (m46) -- ($(m46)!0.9!270:(v6)$) node [sloped,midway,below] {$n_F$} ;
\filldraw [fill=white,draw=red] (m46) circle (3pt); \node[left=-1pt] at (m46) {\color{red}$x_F$};
\node[red] at ($(v4)!0.6!40:(m46)$) {$F$};

\draw[blue,ultra thick] (v8) -- (v9);
\draw[blue,thick,->] (m89) -- ($(m89)!0.7!270:(v9)$) node [sloped,midway,below] {$n_{F'}$} ;
\filldraw [fill=white,draw=blue] (m89) circle (3pt); \node[above left=-3pt] at (m89) {\color{blue}$x_{F'}$};
\node[blue] at ($(v8)!0.25!20:(v9)$) {$F'$};


\node[font=\small] at (3,-1) {b) Unstructured triangular}; 

\end{tikzpicture}
\end{center}
\caption{The Voronoi finite volume method uses cells dual to a Delauney triangulation shown in subfigure a). The discontinuous Galerkin method uses triangular cells defined by a conforming triangular mesh shown in subfigure b).}
\label{fig:fv}
\end{figure}
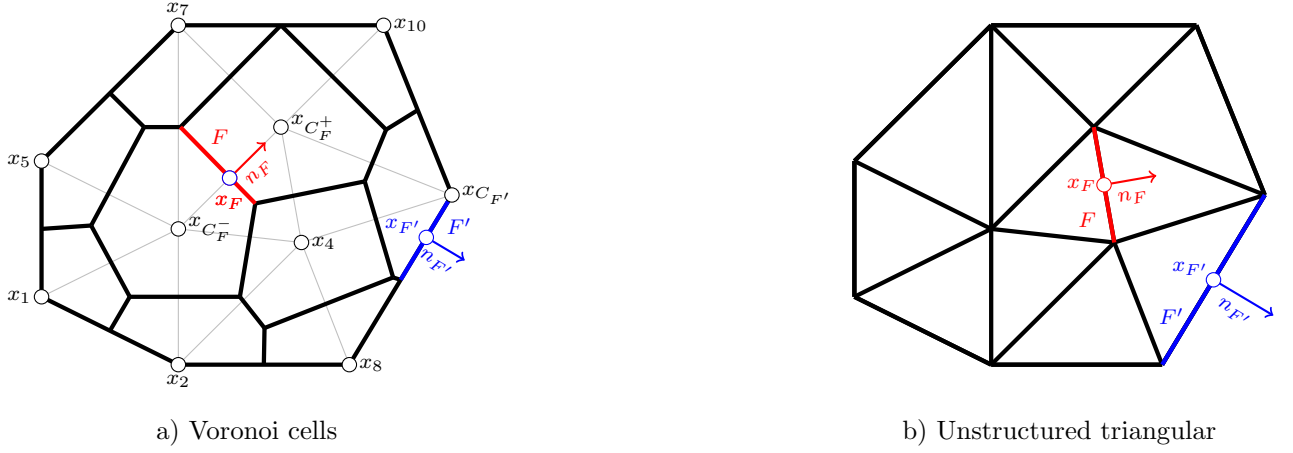

In the Voronoi finite volume method (VFVM) the cells $\mathcal{C}_h$ are given as the Voronoi diagram dual to a Delauney triangulation. For our purposes a Delauney triangulation is a conforming mesh $\mathcal{E}_h = \{E_1, \ldots, E_K\}$ consisting of open triangular elements $E$ partitioning the domain $\Omega$ where the largest angle in any triangle $E\in\mathcal{E}_h$ is not larger than $\pi/2$. The vertices of the mesh $\mathcal{E}_h$ are denoted by $\mathcal{X}_h = \{x_1, \ldots, x_M\}$. Then the Voronoi cell $C_m$ is given by all points in $\Omega$ that are closer to vertex $x_m$ than to any other vertex and $x_m$ is the center of the cell. The intersection $\partial C_i \cap \partial C_j$ of two neighboring cells is a straight line segment formed by the perpendicular bisector of the edge connecting $x_i$ and $x_j$. In Figure~\ref{fig:fv} a) the Delauney triangulation is shown in light gray and the corresponding Voronoi cells are shown in black.   

Finite volume methods approximate the solution using cell-wise constant functions
\begin{equation*}
W^0_h =\left\lbrace  w \in L^2(\Omega): w\vert_C = \text{const} ~\forall C \in \mathcal{C}_h \right\rbrace
\end{equation*}
and we also assume $b\in W_h^0$ for the bathymmetry. Dirichlet boundary conditions are incorporated in strong form by employing the space
\begin{equation*}
W^0_{h,0} =\left\lbrace  w \in W_h^0 : w(x_C) = 0 ~\forall C \in \mathcal{C}_h^D\right\rbrace \subset W_h^0 .
\end{equation*}
The duality of Voronoi diagram and Delauney triangulation allows for a natural postprocessing $\Pi_h : W_h^0 \to V_h^1$ with $$V_h^1 = \{ v\in C^0(\overline\Omega) : v|_{E}\in \mathbb{P}_1 \, \forall E\in\mathcal{E}_h\}$$ the space of piecewise linear conforming finite element functions and $\mathbb{P}_1 = \text{span}\{1,x_1,x_2\}$. Then $\Pi_h$ is the nodal interpolation operator
\begin{equation*}
v = \Pi_h w \quad \Leftrightarrow \quad v(x) = w(x) \, \forall x \in  \mathcal{X}_h .
\end{equation*}

In order to formulate the finite volume method in a concise way, which also highlights the relation to the discontinuous Galerkin method introduced below, we need to introduce some convenient notation. At a point $x \in F$ on an interior edge $F$ a function $w\in W^0_h$ is two-valued:
$w^-(x) = \lim_{\epsilon\to 0+} w(x-\epsilon n_F)$ and $w^+(x) = \lim_{\epsilon\to 0+} w(x+\epsilon n_F)$. The jump of a function  $w\in W^0_h$ at a point $x\in F$ is then defined as
\begin{equation}
\label{eq:jump}
\llbracket w \rrbracket(x)  = w^-(x) - w^+(x)
\end{equation}
and the arithmetic and harmonic average of $w$ at a point $x\in F$ are
\begin{align}
\{ w\}(x) &= \frac{w^-(x) + w^+(x)}{2}, & \langle w \rangle(x)  &= \frac{ 2 w^-(x) w^+(x) }{w^-(x) + w^+(x)} .
\end{align}
For conveniece, jump and averages at $x\in\partial\Omega$ are defined as
\begin{align*}
\llbracket w \rrbracket (x) = \{w\}(x) = \langle w\rangle(x) = w(x).
\end{align*}

\subsubsection{Semi-discrete weak formulation}

The derivation of the VFVM follows a method of lines approach with semi-discretization in space first. For any test function $w\in W^0_{h,0}$ and $t\in\Sigma$, 
the solution of equation \eqref{model_problem_peter} satisfies
\begin{equation} \label{eq:basic_identity}
\begin{split}
\int_\Omega (\partial_t u &+ \nabla \cdot q) w \,dx \\
&= d_t \int_\Omega u w \,dx + \sum_{C\in\mathcal{C}_h}  \int_C (\nabla \cdot q) w \,dx\\
&= d_t \int_\Omega u w \,dx + \sum_{C\in\mathcal{C}_h}  \int_{\partial C} q \cdot n\, w \,ds\\
&= d_t \int_\Omega u w \,dx + \sum_{F\in\mathcal{F}_h^{i}} \int_F q\cdot n_F \, \llbracket w \rrbracket \,ds + \sum_{F\in\mathcal{F}^N_h} \int_F j w \,ds = \int_\Omega f w \,dx .
\end{split}
\end{equation}
Here we used cell-wise Gauß' theorem, the fact that $q\cdot n_F$ is continuous for the exact solution, the definition of the jump and $w=0$ on $\mathcal{F}_h^D$.

The next step is to make the ansatz $u_h(\cdot,t)\in W_h^0$ with $u_h(x_C,t) = g(x_C,t)$ for all $C\in\mathcal{C}_h^D$ and to introduce numerical fluxes $Q_F$ approximating
$q\cdot n_F$. Then the semi-discrete numerical solution reads
\begin{equation}\label{Weak form FV}
\begin{split}
d_t \int_\Omega u_h w \,dx &+ \sum_{F\in\mathcal{F}_h^{i}} \int_F Q_F(t) \, \llbracket w \rrbracket \,ds =\\
&\int_\Omega f w \,dx - \sum_{F\in\mathcal{F}^N_h} \int_F j w \,ds \quad \forall t\in\Sigma, \forall w\in W_{h,0}^0 .
\end{split}
\end{equation}

\subsubsection{Numerical fluxes}

\begin{figure}
\begin{center}
\begin{tikzpicture}[scale=0.9]
\coordinate (v1) at (0,0);
\coordinate (v2) at (5,0);
\coordinate (v3) at (2,3);
\coordinate (v4) at (3,-3);
\coordinate (F1) at (2.5,0.6);
\coordinate (F2) at (2.5,-0.6);
\coordinate (X13) at (1,1.5);
\coordinate (X23) at (3.5,1.5);
\coordinate (X14) at (1.5,-1.5);
\coordinate (X24) at (4,-1.5);
\draw (v1) -- (v2);
\draw (v2) -- (v3);
\draw (v3) -- (v1);
\draw (v1) -- (v4);
\draw (v4) -- (v2);
\draw[ultra thick] (F1) -- (F2);
\draw[ultra thick] (F1) -- (X13);
\draw[ultra thick] (F1) -- (X23);
\draw[ultra thick] (F2) -- (X14);
\draw[ultra thick] (F2) -- (X24);
\node[left] at (2.5,0.1) {$F$};
\node at (2,2.0) {$E_{F,1}$};
\node at (3,-2.0) {$E_{F,2}$};
\filldraw  [fill=white,draw=black]  (v1) circle (3pt);
\node[left] at (v1) {$x_{C_F^-}$};
\filldraw  [fill=white,draw=black]  (v2) circle (3pt);
\node[right] at (v2) {$x_{C_F^+}$};
\filldraw  [fill=white,draw=black]  (v3) circle (3pt);
\filldraw  [fill=white,draw=black]  (v4) circle (3pt);
\node[right] at (2.5,0.3) {$|F\cap E_{F,1}|$};
\node[right] at (2.5,-0.3) {$|F\cap E_{F,2}|$};
\end{tikzpicture}
\end{center}
\caption{Notation used in the definition of the flux $Q_F^{\uparrow,full}$ in Equation \eqref{eq:Q_up_full}.}
\label{fig:full_gradient}
\end{figure}
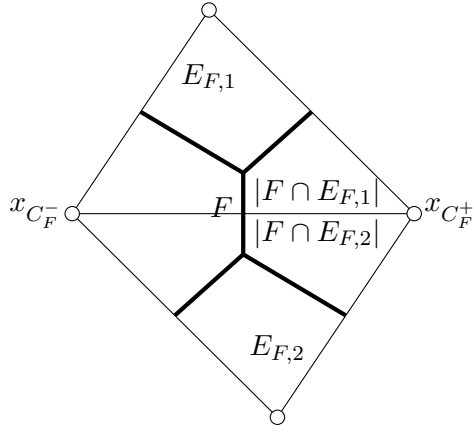

It remains to define the numerical fluxes $Q_F$. Since $F$ is a perpendicular bisector,
\begin{equation*}
\frac{\llbracket u_h(\cdot,t) \rrbracket (x_F)}{\| x_{C^+_F} - x_{C^-_F} \|} = - \nabla u(x_F) \cdot n_F + \mathcal{O}(h^2)
\end{equation*}
is a second order approximation of the directional gradient. 
Following \cite{hunter2005adaptive} we set 
$$ b_F = \max\left(b(x_{C_F^+}),b(x_{C_F^-})\right) $$
and introduce the non-negative upwind water height as
\begin{align} \label{eq:upwind_water_height}
H_F^\uparrow(t) = \max\left( 
\max\left( u_h(x_{C_F^+},t),u_h(x_{C_F^-},t) \right)- b_F,0 \right) .
\end{align}
Such a construction also appears as the interface hydrostatic reconstruction in well-balanced finite volume methods for the full shallow water equations \cite{AUDUSSE2005311}.
For the conductivity we choose the harmonic average
\begin{align*}
K_F = \langle K \rangle (x_F)
\end{align*}
and for the gradient in the friction law we employ the piecewise linear reconstruction on the Delauney triangulation. 
Every $F\in\mathcal{F}_h^i$ is associated with two triangles $E_{F,1}, E_{F,2}$ which are adjacent to the edge connecting $x_{C_F^-}$ and $x_{C_F^+}$, see Figure \ref{fig:full_gradient}. On each of these triangles we can approximate the norm of the gradient as
\begin{align*}
G_{F,1} &= \| \nabla (\Pi_h u_h) |_{E_{F,1}}\| + \epsilon, &
G_{F,2} &= \| \nabla (\Pi_h u_h) |_{E_{F,2}}\| + \epsilon,
\end{align*}
using piecewise linear interpolation of the values in the vertices in each of the two triangles. Denoting by $|F\cap E_{F,1}|$, $|F\cap E_{F,2}|$ the lengths of the parts of $F$ lying in the respective triangle and by $|F|$ the total length of $F$ we define the four-point numerical flux
\begin{align}
Q_F &= K_F \left(H_F^\uparrow\right)^\alpha
\left( \frac{|F\cap E_{F,1}|}{(G_{F,1})^{1-\gamma}|F|} + \frac{|F\cap E_{F,2}|}{(G_{F,2})^{1-\gamma}|F|}  \right)
\frac{\llbracket u_h(\cdot,t) \rrbracket (x_F)}{\| x_{C^+_F} - x_{C^-_F} \|} . \label{eq:Q_up_full}
\end{align}
In Theorem \ref{thm:nonnegativity_VFVM} below we show that this flux leads to a provably nonnegative discrete solution.

\subsubsection{Fully-discrete formulation}

For simplicity we proceed with the implicit Euler method for time discretization to formulate the fully discrete scheme and observe that all functions to be integrated are piecewise constant on cells or triangles respectively. Then, given $u_h^{n-1}\in W_h^0$, $n>0$, we seek $u_h^{n}\in W_h^0$ with $u_h^n(x_C) = g(x_C,t^n)$ for all $C\in\mathcal{C}_h^D$ such that
\begin{equation}\label{fully_discrete_fv}
\begin{split}
\sum_{C \in \mathcal{C}^{iN}} & u_h^n(x_C) w(x_C) |C| 
+ \Delta t^n \sum_{F\in\mathcal{F}_h^{i}} Q_F(t^n) \, \llbracket w \rrbracket(x_F) |F| = \sum_{C \in \mathcal{C}^{iN}} u_h^{n-1}(x_C) w(x_C) |C| \\
&
+ \Delta t^n \sum_{C \in \mathcal{C}^{iN}} f(x_C) w(x_C) |C| 
- \Delta t^n \sum_{F\in\mathcal{F}^N_h} j(x_F,t^n) w(x_F) |F| \quad \forall w\in W_{h,0}^0 .
\end{split}
\end{equation}
Other implicit time-stepping schemes, such as diagonally-implicit Runge-Kutta methods from \cite{alexander1977diagonally} can be used in a similar way.

\begin{theorem}[Nonnegativity of water height in the VFVM] \label{thm:nonnegativity_VFVM}
    Consider the fully discrete scheme \eqref{fully_discrete_fv} with homogeneous Neumann boundary conditions $j(x,t)=0$, no source or sink term $f(x,t) = 0$, and nonnegative initial condition $u^{0}_h(x_C) \geq b(x_C)$ $\forall C\in\mathcal{C}_h$. Then the fully discrete solution, if it exists, satisfies
\begin{equation*}
u^n_h(x_C) \geq b(x_C), \qquad \forall n>0, \ \forall C \in \mathcal{C}_h^{iN}.
\end{equation*}
\begin{proof}
The proof is by induction over the time steps and using a contradiction argument. The initial condition has nonnegative water height by assumption. For the induction step assume now that water heights in $u_h^{n-1}$ are nonnegative. 
Suppose now that that 
\begin{equation}
u^n_h(x_C) - b(x_C) < 0   \label{contradiction_assumption}
\end{equation}
for some $C\in\mathcal{C}_h$.
Inserting the test function $w = \chi_C$ (being one on cell $C$ and zero else) and assuming without loss of generality that all unit normal vectors $n_F$ are oriented such that they point outward of $C$ (then $\llbracket w \rrbracket = 1$) \eqref{fully_discrete_fv} reduces to
\begin{align*}
u^n_h(x_C) = u^{n-1}_h(x_C) - \frac{\Delta t^n}{\vert C \vert} \sum_{F\in\mathcal{F}_h^{i} \cap \partial C} Q_{F}(t^n) \vert F \vert .
\end{align*}
Subtracting $b(x_C)$ in both sides of the equation above yields
\begin{align*}
0 > u^n_h(x_C) - b(x_C) = u^{n-1}_h(x_C) - b(x_C) - \frac{\Delta t^n}{\vert C \vert} \sum_{F\in\mathcal{F}_h^{i} \cap \partial C} Q_{F}(t^n) \vert F \vert
\end{align*}
which is equivalent to
\begin{align*}
    0 \leq u^{n-1}_h(x_C) - b(x_C) < \frac{\Delta t^n}{\vert C \vert} \sum_{F\in\mathcal{F}_h^{i}\cap \partial C} Q_{F}(t^n) \vert F \vert
\end{align*}
since $u^{n-1}_h(x_C) - b(x_C) \geq 0$ by induction assumption. Now $\vert C \vert, \vert F \vert, \Delta t^n > 0$ and this implies there must be at least one $F\in\mathcal{F}_h^{i}\cap \partial C$ with $Q_{F}(t^n)>0$, i.e. $F$ is an outflow edge of $C$ since $n_F$ points outward of $C$. The numerical flux \eqref{eq:Q_up_full} then implies that $H_F^\uparrow>0$ and $u_h^n(x_{C}) > u_h^n(x_{C_F^+})$ (note that $x_{C} = x_{C_F^-}$). The definition of the upwind water height \eqref{eq:upwind_water_height} then implies $$0 < H_F^\uparrow = u_h^n(x_C) - \max(b(x_C),b(x_{C_F^+})).$$
Now there are two cases
\begin{itemize}
\item[i)] $b(x_C)\geq b(x_{C_F^+})$, then we have $u_h^n(x_C)-b(x_C)>0$ which contradicts assumption \eqref{contradiction_assumption}.
\item[ii)] $b(x_C) < b(x_{C_F^+})$, but then $0 < u_h^n(x_C) - \max(b(x_C),b(x_{C_F^+})) = u_h^n(x_C) - b(x_{C_F^+}) < u_h^n(x_C) - b(x_C)$ which also contradicts assumption \eqref{contradiction_assumption}.
\end{itemize}
Since $C$ was an arbitrary cell with negative water height we have shown by contradiction that $u^n_h(x_C) \geq b(x_C)$ for all $C\in\mathcal{C}_h^{iN}$ and $n>0$.
\end{proof}
\end{theorem}

The proof shows that the upwind water height and the two-point evaluation of the gradient in normal direction are both essential to achieve nonnegativity.

\subsection{Discontinuous Galerkin schemes}

The VFVM with upwind flux is simple to implement and yields provably nonnegative water height. Nevertheless, improvement is needed in several respects:
\begin{enumerate}[a)]
\item The possibility to use continuous and discontinuous bathymmetry would be advantageous to handle overland flows and urban flooding with the same scheme.
\item The Delauney triangulation is rather restrictive. In particular, it is not amenable to adaptive local mesh refinement except in the case of equilateral, right-angled triangles with newest vertex bisection.
\item The upwind water height leads to a first-order accurate scheme when second-order accuracy would be desirable.
\end{enumerate}
We now introduce a DG method that improves on all three aspects. First we develop a basic DG method which is then extended by a cut-cell approach improving the behavior at the wet/dry front.

\subsubsection{Semi-discrete weak formulation}

In the DG method the set of cells $\mathcal{C}_h$ coincides now with the elements of the triangular mesh $\mathcal{E}_h$, so $M=K$, see Figure \ref{fig:fv} part b. Moreover, the mesh is not required to be Delaunay but for simplicity we assume that it is conforming (but this condition could be relaxed as well). Now the approximation space is the space of element-wise linear functions:
\begin{equation*}
W^1_h =\left\lbrace  w \in L^2(\Omega): w\vert_C \in \mathbb{P}_1 ~\forall C \in \mathcal{C}_h\right\rbrace 
\end{equation*}
which is used for the solution, test functions and bathymmetry. Note that $V_h^1 \subset W_h^1$ as well as $W_h^0 \subset W_h^1$ where $W_h^0$ is now defined w.r.t. the triangular cells.

Multiplying Equation \eqref{eq:conservation} by a test function $w\in W_h^1$, integrating over the domain and using element-wise integration by parts one arrives at
\begin{equation} \label{eq:basic_identity_2}
\begin{split}
d_t \int_\Omega u w \,dx &- \sum_{C\in\mathcal{C}_h}  \int_C q \cdot  \nabla w \,dx + \sum_{F\in\mathcal{F}_h^{iD}} \int_F q\cdot n_F \llbracket w \rrbracket \,ds \\ 
&= \int_\Omega f w \,dx - \sum_{F\in\mathcal{F}^N_h} \int_F j w \,ds \qquad \forall t\in\Sigma, \forall w \in W_h^1.
\end{split}
\end{equation}
In contrast to the VFVM, Dirichlet boundary conditions are incorporated in a weak form and the test function $w$ is not zero at the Dirichlet boundary.

The semidiscrete weak formulation for the numerical solution $u_h(\cdot,t)\in W_h^1$ is then obtained by introducing an appropriate numerical flux $Q_F$ and adding stabilization terms $P_F$, $M_F$ \cite{riviere2008discontinuous}:
\begin{equation} \label{eq:basic_identity_3}
\begin{split}
d_t \int_\Omega u_h & w \,dx - \sum_{C\in\mathcal{C}_h}  \int_C q \cdot  \nabla w \,dx + \sum_{F\in\mathcal{F}_h^{iD}} \int_F Q_F(u_h,t) \llbracket w \rrbracket \,ds \\
&\qquad + \sum_{F\in\mathcal{F}_h^{iD}} \int_F P_F(w,t) \llbracket u_h \rrbracket \,ds
+ \sum_{F\in\mathcal{F}_h^{iD}} \int_F M_F(u_h,t) \llbracket u_h \rrbracket \llbracket w \rrbracket\,ds \\
& = \int_\Omega f w \,dx - \sum_{F\in\mathcal{F}^N_h} \int_F j w \,ds \\
&\qquad+ \sum_{F\in\mathcal{F}_h^{D}} \int_F M_F(u_h,t) g w \,ds 
\quad \forall t\in\Sigma, \forall w \in W_h^1.
\end{split}
\end{equation}

\subsubsection{Numerical flux and stabilization} \label{sec:UDGM} \label{sec:analysis}

For any edge $F\in\mathcal{F}_h^{iD}$ we first define the flux direction as
\begin{equation}\label{eq:flow_direction}
D_F(x,t) = \begin{cases}
- \langle K \rangle \{  \nabla u_h \cdot n_F \}(x,t) + \frac{\sigma_F}{h_F} \langle K \rangle  \llbracket u_h \rrbracket(x,t) & F \in \mathcal{F}_h^i \\
- K \nabla u_h(x,t) \cdot n_F + \frac{\sigma_F}{h_F} K (u_h(x,t)-g(x,t))  & F \in \mathcal{F}_h^D
\end{cases},
\end{equation}
where $\sigma_F$ is a user chosen penalty factor and $h_F$ is the length of $F$.
Then, depending on the flux direction, the upwind water height is determined:
\begin{align*} \label{classical flux formulation DG}
H_F^\uparrow(x,t) &= \begin{cases}
\max(u_h^-(x,t) - b_F(x),0), & D_F(x,t) \geq 0 \\
\max(u_h^+(x,t) - b_F(x),0) & D_F(x,t) < 0,
\end{cases} &&F\in \mathcal{F}_h^i, \\ 
H_F^\uparrow(x,t) &= \begin{cases}
\max(u_h^-(x,t) - b(x),0), & D_F(x,t) \geq 0 \\
\max(g(x,t) - b(x),0) & D_F(x,t) < 0
\end{cases} &&F\in \mathcal{F}_h^D.
\end{align*}
where on interior edges $b_F(x) = \max(b^-(x),b^+(x))$ is the larger of the two values of bathymmetry at an edge. The numerical gradient is evaluated as the average:
\begin{equation*}
G_F(x,t) = \begin{cases} 
\left\|\{ \nabla u_h(x,t) \} \right\| + \epsilon & F \in \mathcal{F}_h^i \\
\left\| \nabla u_h(x,t) \right\| + \epsilon & F \in \mathcal{F}_h^D
\end{cases} .
\end{equation*}
With this in place the numerical flux and stabilization terms are defined as
\begin{align}
Q_F(u_h,t) &= - \frac{H_F^\uparrow(x,t)^\alpha}{G_F(x,t)^{1-\gamma}} \langle K \rangle \{\nabla u_h(x,t) \cdot n_F \}, \label{QF_DG}\\
P_F(w,t) &= - \frac{H_F^\uparrow(x,t)^\alpha}{G_F(x,t)^{1-\gamma}} \langle K \rangle \{\nabla w(x,t) \cdot n_F \}, \label{PF_DG}\\
M_F(u_h,t) &= \frac{H_F^\uparrow(x,t)^\alpha}{G_F(x,t)^{1-\gamma}} \frac{\sigma_F}{h_F} \langle K \rangle. \label{MF_DG}
\end{align}
Observe that upon inserting a test function $w\in W_h^0$ this formulation results in a total flux $\frac{H_F^\uparrow(x,t)^\alpha}{G_F(x,t)^{1-\gamma}} D_F(x,t)$ on the edge $F$.
The scheme is an extension of the symmetric interior penalty discontinuous Galerkin method introduced in \cite{ern2009discontinuous} and used, e.g. in \cite{bastian2014fully} for nonlinear parabolic equations.

A fully discrete formulation is obtained by inserting \eqref{QF_DG} -- \eqref{MF_DG} into  \eqref{eq:basic_identity_3}, employing an appropriate implicit time integration and use of numerical quadrature for evaluating the integrals.

\subsubsection{Cut-cell formulation} \label{sec:CUDGM}

In this subsection we use a cut-cell DG (or unfitted DG) approach to ensure non-negativity of the solution near the wet-dry front. Cut-cell (or unfitted) finite element methods are an established approach for the numerical solution of partial differential equations in complex domains which are not resolved by the mesh, see e.g. \cite{HANSBO20025537,Bastian2009,Heimann2013,BURMAN20101217,Burman2015,BADIA2018533,engwersisc2020}. While often the complex domain is given externally, e.g. through a level set function \cite{Bastian2009,Heimann2013,Burman2015}, the complex domain here is the time-dependent support of the solution $\Omega^+(t) = \supp u(t)$.

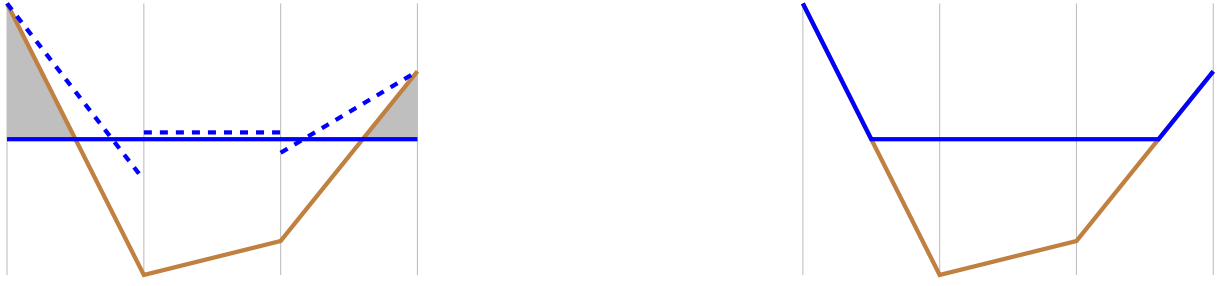
\begin{figure}
\begin{center}
\begin{tikzpicture}[scale=0.9]
\draw[lightgray,thin] (0,0) -- (0,4);
\draw[lightgray,thin] (2,0) -- (2,4);
\draw[lightgray,thin] (4,0) -- (4,4);
\draw[lightgray,thin] (6,0) -- (6,4);
\coordinate (b1) at (0,4);
\coordinate (b2) at (2,0);
\coordinate (b3) at (4,0.5);
\coordinate (b4) at (6,3);
\coordinate (u1) at (0,2);
\coordinate (u2) at (2,2);
\coordinate (u2x) at (1,2);
\coordinate (u3) at (4,2);
\coordinate (u3x) at (5.2,2);
\coordinate (u4) at (6,2);
\fill[lightgray] (0,2) -- (u2x) -- (b1) -- cycle;
\fill[lightgray] (u3x) -- (6,2) -- (b4) -- cycle;
\draw[ultra thick,draw=brown] (b1) -- (b2) -- (b3) -- (b4);
\draw[ultra thick,draw=blue] (u1) -- (u2) -- (u3) -- (u4);
\draw[ultra thick,draw=blue,dashed] (0,4) -- (2,1.4); 
\draw[ultra thick,draw=blue,dashed] (2,2.1) -- (4,2.1); 
\draw[ultra thick,draw=blue,dashed] (4,1.8) -- (6,3.0); 
\node at (3,-1) {a) water surface in a linear bathymmetry};
\end{tikzpicture}
\hfill
\begin{tikzpicture}[scale=0.9]
\draw[lightgray,thin] (0,0) -- (0,4);
\draw[lightgray,thin] (2,0) -- (2,4);
\draw[lightgray,thin] (4,0) -- (4,4);
\draw[lightgray,thin] (6,0) -- (6,4);
\coordinate (b1) at (0,4);
\coordinate (b2) at (2,0);
\coordinate (b3) at (4,0.5);
\coordinate (b4) at (6,3);
\draw[ultra thick,draw=brown] (b1) -- (b2) -- (b3) -- (b4);
\coordinate (u1) at (0,4);
\coordinate (u2x) at (1,2);
\coordinate (u2) at (2,2);
\coordinate (u3) at (4,2);
\coordinate (u3x) at (5.2,2);
\coordinate (u4) at (6,3);
\draw[ultra thick,draw=blue] (u1) -- (u2x) -- (u2) -- (u3) -- (u3x) -- (u4);
\node at (3,-1) {b) desired water surface};
\end{tikzpicture}
\end{center}
\caption{Motivation for cut-cell DG.}
\label{fig:cut_cell_motivation}
\end{figure}

Piecewise linear bathymmetry and solution may pose a problem with respect to the nonnegativity of the water height. We illustrate this in Figure \ref{fig:cut_cell_motivation} for a one-dimensional situation. Subfigure a) shows three elements with a continuous linear bathymmetry in solid brown. A constant, stationary water level (lake at rest) is shown in solid blue. It is fine to represent the lake surface but exhibits negative water height in the gray shaded regions. Moreover, this ``negative mass'' will be compensated by excess positive mass elsewhere since the DG scheme is locally mass conservative. Insisting on a nonnegative water height, e.g. through the use of limiters or a variational inequality formulation, might lead to the dashed blue solution which now has non-negative water height but is not a good solution at the wet/dry front. The desired solution, which is not a function in $W_h^1$, is the function shown in subfigure b). This solution is obtained by the cut-cell DG method described in this section.

In order to formulate the cut-cell DG method, we define the nonlinear space of admissible water heights
\begin{align} \label{nnp-set}
     A_{h}^1 = \{u \in L_2(\Omega)\, \vert \, u(x) = \max(0,v(x)),\, v\in W^1_h \},
\end{align}
and seek the fully discrete solution as $u_h(t) = b+a_h(t)$ with $a_h(t)\in A_h^1$. Note that the solution shown in Figure \ref{fig:cut_cell_motivation} b) is exactly of this form.

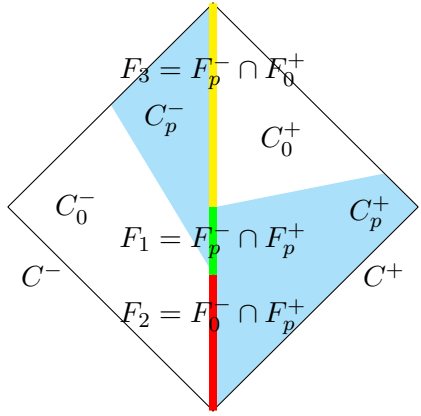
\begin{figure}[t]
    \centering
    \begin{tikzpicture}[scale = 0.9]
    \draw[fill=cyan!30,draw=white] (0,-1) -- (-1.5,1.5) -- (0,3) -- (0,-1);    
    \draw[fill=cyan!30,draw=white] (0,0.0) -- (0,-3) -- (3,0) -- (2.5,0.5) --(0,0.0);
    \draw (-3,0) -- (0,3) -- (0,-3) -- (-3,0) ;
    \draw (3,0) -- (0,3) -- (0,-3) -- (3,0) ;

    \node at (-2.5,-1) {$C^-$};
    \node at (2.5,-1) {$C^+$};
    \node at (-0.7,1.3) {$C^-_p$};
    \node at (2.3,-0.1) {$C^+_p$};
    \node at (-2,0) {$C^-_0$};
    \node at (1,1) {$C^+_0$};



    \draw[line width=3pt,red] (0,-3) -- (0,-1);
    \draw[line width=3pt,green] (0,-1) -- (0,0);
    \draw[line width=3pt,yellow] (0,0) -- (0,3);

    \node at (0,-0.5) {$F_1 = F_p^- \cap F_p^+$};
    \node at (0,-1.6) {$F_2 = F_0^- \cap F_p^+$};
    \node at (0,2) {$F_3 = F_p^- \cap F_0^+$};

    
    \end{tikzpicture}
    \caption{Illustration of cut cells and cut edges.}
    \label{fig:cut edge}
\end{figure}

As a consequence, the storage term in the semi-discrete weak formulation now reads for any test function $w\in W_h^1$:
\begin{equation*}
\begin{split}
    d_t \int_\Omega u_h(t) w \,dx &= 
    d_t \int_\Omega (b + a_h(t)) w \,dx\\
    &= \sum_{C\in C_h} d_t \int_C a_h(t) w \,dx =
    \sum_{C\in C_h} d_t \int_{C_p(t)} a_h(t) w \,dx
\end{split}
\end{equation*}
where $C_p(t) = \supp{a_h(t)} \cap C$ is the time dependent support of the function $a_h(t)$ in the triangle $C$. If $C_p(t) \subset C$ we call $C_p(t)$ a cut-cell. The integral over $C_p(t)$ can be carried out exactly, because, for a fixed $t$, $a_h$ is linear on $C_p$ and $C_p$ is either of triangular or quadrilateral form or empty. Figure \ref{fig:cut edge} shows two cells with their cut cells as the blue shaded regions. The second volume term in \eqref{eq:basic_identity_3} is treated in the same way, since $q$ has also the support $C_p$ on cell $C$. 

It remains to explain the terms involving integrals over edges. For these observe that $u_h(t)=b + a_h(t)$ splits each cell $C\in\mathcal{C}_h$ into two parts $C = C_p(t) \cup C_0(t)$, where $C_0 = C\setminus C_p(t)$, as indicated in Figure \ref{fig:cut edge} (we omitted the time dependence for brevity). Now consider an interior edge $F\in\mathcal{F}_h^i$ with its two adjacent cells $C_F^-, C_F^+$. The splitting into cut cells induces the following splitting of the edge $F$:
\begin{align*}
F_p^- &= F \cap \partial C_{F,p}^-, & F_0^- &= F \cap \partial C_{F,0}^-, &
F_p^+ &= F \cap \partial C_{F,p}^+, & F_0^+ &= F \cap \partial C_{F,0}^+.
\end{align*}
In order to carry out the edge integrals as accurately as possible we partition the edge $F$ into the following cut edges:
\begin{enumerate}[1)]
    \item $F_1 = F_p^- \cap F_p^+$ shown in green in Figure \ref{fig:cut edge}: on both sides the water height is positive and the edge integral can be evaluated as before.
    \item $F_2 = F_0^- \cap F_p^+$ shown in red in Figure \ref{fig:cut edge}: here only flux from $C^+$ to $C^-$ can occur.
    \item $F_3 = F_p^- \cap F_0^+$ shown in yellow in Figure \ref{fig:cut edge}: here only flux from $C^-$ to $C^+$ can occur.
    \item On $F_0^- \cap F_0^+$ the flux is zero because water height is zero on both sides of the edge.
\end{enumerate}
For later use we set $$F_p(t) = F_1(t) \cup F_2(t) \cup F_3(t)$$
where we have made explicit again the time dependence of these domains.
Within each cut edge $F_\alpha$, $\alpha\in\{1,2,3\}$, the flux direction $D_F$ \eqref{eq:flow_direction} needs to be determined. Since $\{\nabla u_h\cdot n_F\}$ is constant for piecewise linears and $\llbracket u_h \rrbracket$ is linear for fixed $\alpha$ there is at most one change of flux direction within each $F_\alpha$. Thus, each edge is partitioned into at most six parts for numerical integration.

The discrete weak form for the cut-cell DG method now reads as follows. Find $u_h(t) \in b + A_h^1$ such that
\begin{equation} \label{eq:cut_cell_DG}
\begin{split}
\sum_{C\in\mathcal{C}_h} & d_t \int_{C_p(t)} u_h w \,dx 
- \sum_{C\in\mathcal{C}_h}  \int_{C_p(t)} q \cdot  \nabla w \,dx + \sum_{F\in\mathcal{F}_h^{iD}} \int_{F_p(t)} Q_F(u_h,t) \llbracket w \rrbracket \,ds \\
&\quad+ \sum_{F\in\mathcal{F}_h^{iD}} \int_{F_p(t)} P_F(w,t) \llbracket u_h \rrbracket \,ds
+ \sum_{F\in\mathcal{F}_h^{iD}} \int_{F_p(t)} M_F(u_h,t) \llbracket u_h \rrbracket \llbracket w \rrbracket\,ds \\
& = \int_\Omega f w \,dx - \sum_{F\in\mathcal{F}^N_h} \int_F j w \,ds \\
&\quad+ \sum_{F\in\mathcal{F}_h^{D}} \int_{F_p(t)} M_F(u_h,t) g w \,ds 
\qquad \forall t\in\Sigma, \forall w \in W_h^1.
\end{split}
\end{equation}

The fully discrete scheme is obtained by chosing an implicit time discretization, e.g. implicit Euler, and solving a nonlinear algebraic system for the non-negative water height $a_h^n$ at time $t^n$. Choosing a DG basis with element-local support, we represent the discrete solution as
$$a_h^n(x) = \sum_{C\in\mathcal{C}_h} \max(0,z_{C,0}\phi_{C,0}(x)+z_{C,1}\phi_{C,1}(x)+z_{C,2}\phi_{C,2}(x)).$$
Clearly, if in a cell the local polynomial becomes negative on the whole cell, the Jacobian matrix of the nonlinear algebraic system will become singular. It will become ill-conditioned if $\supp a_h^n \cap C$ becomes very small. This is a manifestation of the ``small cell problem'' well-known in cut-cell methods and several approaches have been devised to solve it. For elliptic and parabolic problems solved implicitly, the ghost penalty approach \cite{BURMAN20101217,Burman2015} adds an additional penalty term near the boundary and the aggregation approach \cite{Heimann2013,BADIA2018533} adds small cells to neighboring cells. A special penalty term for explicit time discretization of hyperbolic PDEs was suggested in \cite{engwersisc2020}. However, \cite{Heimann2013} reports that aggregating and not aggregating did not make a huge difference when solving a two-phase Navier-Stokes problem.

Here we take a regularization approach ensuring that a cell never runs completely dry. The idea of this regularization is to replace $H_F^\uparrow$ in \eqref{QF_DG} -- \eqref{MF_DG} by its regularization $\nu(H_F^\uparrow)$
reducing the flux to zero before $H_F^\uparrow$ becomes zero. 
Specifically, we set for two regularization parameters $0 < \delta_1 < \delta_2$
\begin{equation} \label{zeroing function}
\nu(H) = 
    \begin{cases}
        0, &H < \delta_1, \\
        p(H), &\delta_1 \leq H < \delta_2, \\
        H, &H \geq \delta_2.
    \end{cases}
\end{equation}
where $p(H)$ is a polynomial of degree three given by 
\begin{equation*}
    p(H) = \frac{(H - \delta_1)^2}{(\delta_2 - \delta_1)^2}\left( (\delta_1 + 2\delta_2) - (\delta_1 + \delta_2)\frac{H - \delta_1}{\delta_2 - \delta_1}\right)
\end{equation*}
ensuring that $p$ is continuously differentiable.
The flux regularization is combined with a minimum water depth $0<\eta_0<\delta_1$ as initial value for dry cells. 
In the convergence tests below the regularization parameters are chosen depending on the mesh size.

\section{Numerical results} \label{result section}

This section presents three numerical examples designed to validate and compare the performance of the three proposed schemes: the Finite Volume Method on Voronoi meshes (VFVM) from Section \ref{sec:VFVM}, the upwind discontinous Galerkin method (UDGM) given in Section \ref{sec:UDGM} and the cut-cell upwind discontinous Galerkin method (CUDGM) described in Section \ref{sec:CUDGM}. 

In Section \ref{sec:barenblatt} we conduct a convergence study for a problem with analytical solution to analyze the accuracy of the proposed schemes including the wet/dry front. In Section \ref{sec:obstacle} we simulate the full DWE model on an inclined plane with an obstacle. This example shows the versatility and accuracy of the cut-cell DG scheme handling the full nonlinearity on continuous and discontinuous bathymmetry.
Finally, in Section \ref{Sec:dam_break} we demonstrate the robustness of the cut-cell DG method for the full non-linearity and a real-world bathymetry. The results show that the cut-cell DG scheme is highly accurate, consistently outperforming the VFVM method. Note that all results of VFVM presented in this section are using $\mathbb{P}_1$ post-processing.

\subsection{Barenblatt example} \label{sec:barenblatt}

We consider a specific case of our model \eqref{model_problem_peter} by setting $\gamma = \alpha = 1$, and $K = 2$. This simplifies the governing equation to:
\begin{align}
    \partial_t u - 2\nabla \cdot ((u-b) \nabla u) = 0 \label{simplified model equation_u}
\end{align}
Exploiting that bathymmetry $b(x)$ does not depend on time, an equivalent formulation in terms of the water depth $H(x,t)$ is:
\begin{align}
    \partial_t H - 2\nabla \cdot (H \nabla (H+b)) = 0 \label{simplified model equation_H}
\end{align}
Now we consider two cases for the bathymmetry $b(x)$. In the first case $b(x) = 0$, the model reduces to the porous medium equation (now $u = H$):
\begin{align}
    \partial_t u - 2\nabla \cdot (u \nabla u) = 0 \,. \label{porous medium equation}
\end{align}
This equation has the well-known Barenblatt analytical solution:
\begin{align}
    u(x,t) = \max \bigg[ 0, t^{-1/2} \bigg( M - \frac{1}{16} \frac{\|x \|^2}{t^{1/2}} \bigg) \bigg] \label{barenblatt m=2}
\end{align}
where we set $M=0.2$ (in general $M$ depends on the mass of the solution). 

In the second case we consider a uniformly inclined plane, where $b(x) = \vec{v} \cdot x$ for a constant vector $\vec{v} \in \mathbb{R}^2$. The $H$-based model then becomes:
\begin{align}
    \partial_t H - 2\nabla \cdot (H \nabla H) - 2\vec{v} \cdot \nabla H = 0 \label{porous medium equation H-form in slope}
\end{align}
adding a linear advection term to the porous medium equation.
The solution to equation~\eqref{porous medium equation H-form in slope} can be derived using the following lemma, which relates it back to the solution of the standard porous medium equation.

\begin{lemma} \label{Barenblatt lemma}
Let $u(x,t)$ be the solution of the porous medium equation
\[ \partial_t u - 2 \nabla \cdot ( u \nabla u) = 0 \,. \]
Then $\hat w(\hat x, \hat t) = u(\hat x + 2 \vec{v} \hat t, \hat t)$ solves
\[ \partial_{\hat t} \hat w - 2 \hat\nabla \cdot (\hat w \hat\nabla \hat w) - 2 \vec{v} \cdot \hat\nabla \hat w = 0 \,,\]
where $\hat\nabla$ means differentiation w.r.t. $\hat x$.
\begin{proof}
Use the transformation of variables $x = \hat x + 2 \vec{v}\hat t$, $t = \hat t$ and the chain rule. \end{proof}
\end{lemma}


\noindent Thus equation \eqref{porous medium equation H-form in slope}
has the analytical solution 
\begin{align}
    H(x ,t) = \max \bigg[ 0, t^{-1/2} \bigg( M - \frac{1}{16} \frac{\| x + 2\vec{v} t \|^2}{t^{1/2}} \bigg) \bigg] . \label{analytical solution}
\end{align}

\begin{figure}
     \centering
     \begin{subfigure}[b]{0.45\textwidth}
         \centering
         \includegraphics[width=\textwidth]{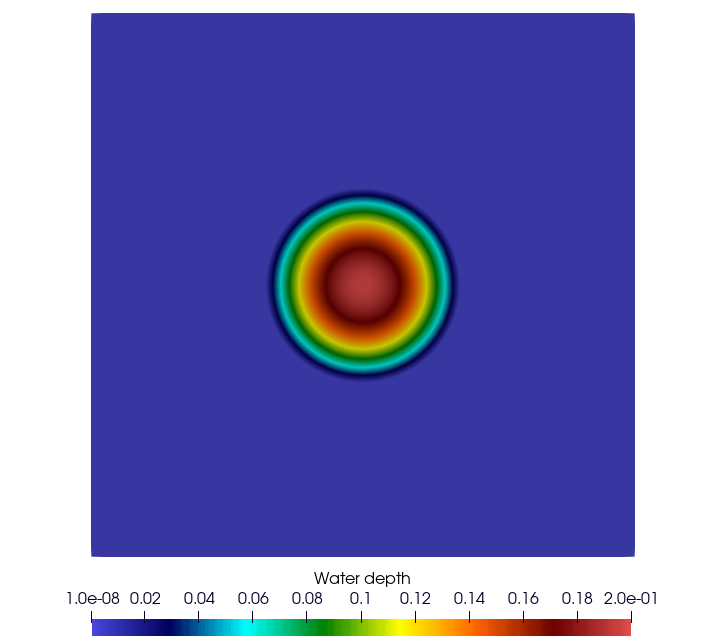}
         \caption{$\vec{v} = (0,0)^T,\, t = 1$}
         \label{exact 2D V00 t0}
     \end{subfigure}
     \hfill
     \begin{subfigure}[b]{0.45\textwidth}
         \centering
         \includegraphics[width=\textwidth]{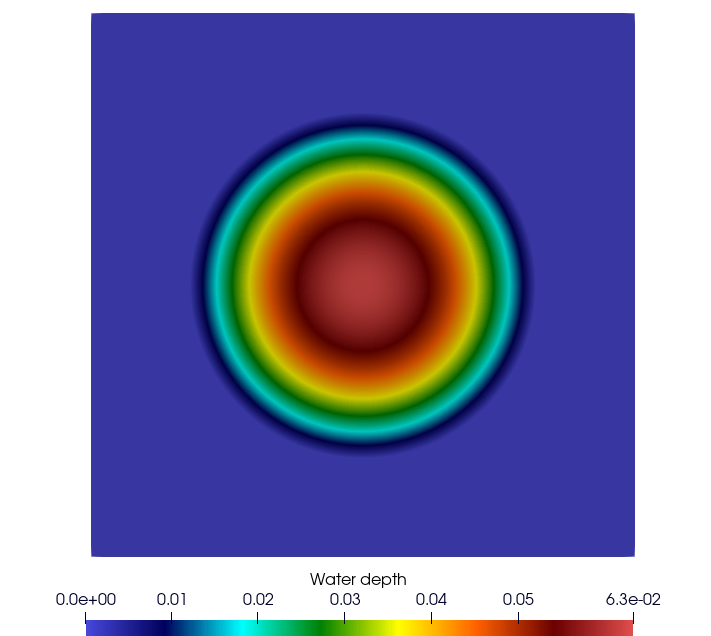}
         \caption{$\vec{v} = (0,0)^T,\, t = 10$}
         \label{exact 2D V00 t1.5}
     \end{subfigure}
     \begin{subfigure}[b]{0.45\textwidth}
         \centering
         \includegraphics[width=\textwidth]{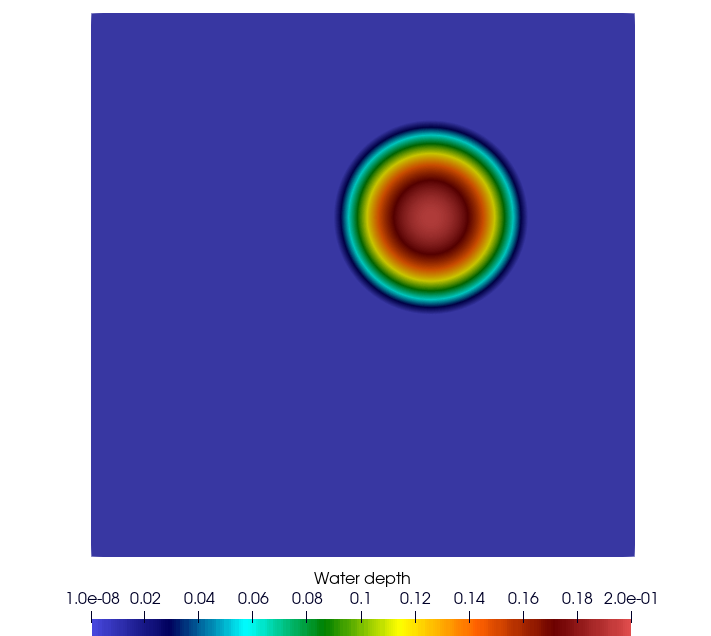}
         \caption{$\vec{v} = (1/2,1/2)^T,\, t = 1$}
         \label{exact 2D V55 t0}
     \end{subfigure}
     \hfill
     \begin{subfigure}[b]{0.45\textwidth}
         \centering
         \includegraphics[width=\textwidth]{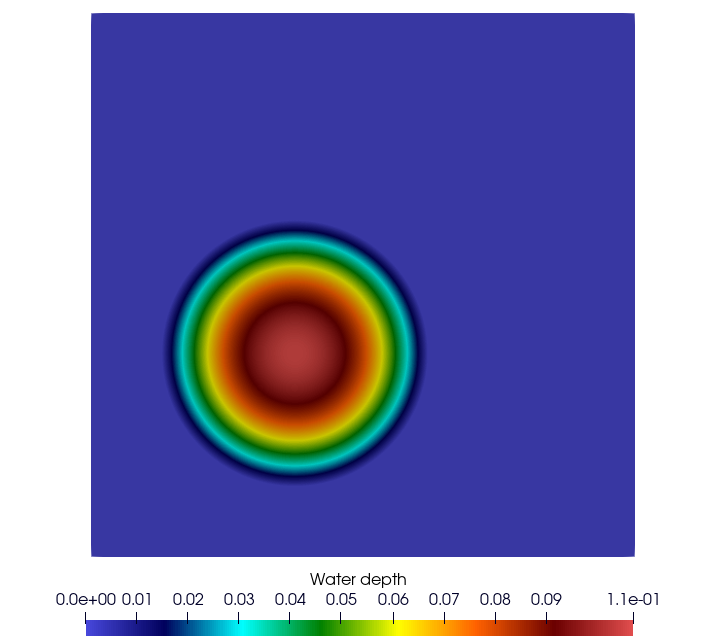}
         \caption{$\vec{v} = (1/2,1/2)^T,\, t = 3.5$}
         \label{exact 2D V55 t1.5}
     \end{subfigure}
    \caption{The analytical solution \eqref{analytical solution} for $\vec{v} = 0$ (first row, corresponding to $b=0$) and $\vec{v} = (1/2,1/2)^T$ (second row, corresponds to $b = \vec{v}\cdot x$).}
    \label{fig:exact 2D}
\end{figure}

\begin{figure}
     \centering
     \begin{subfigure}[b]{0.4\textwidth}
         \centering
         \includegraphics[width=\textwidth]{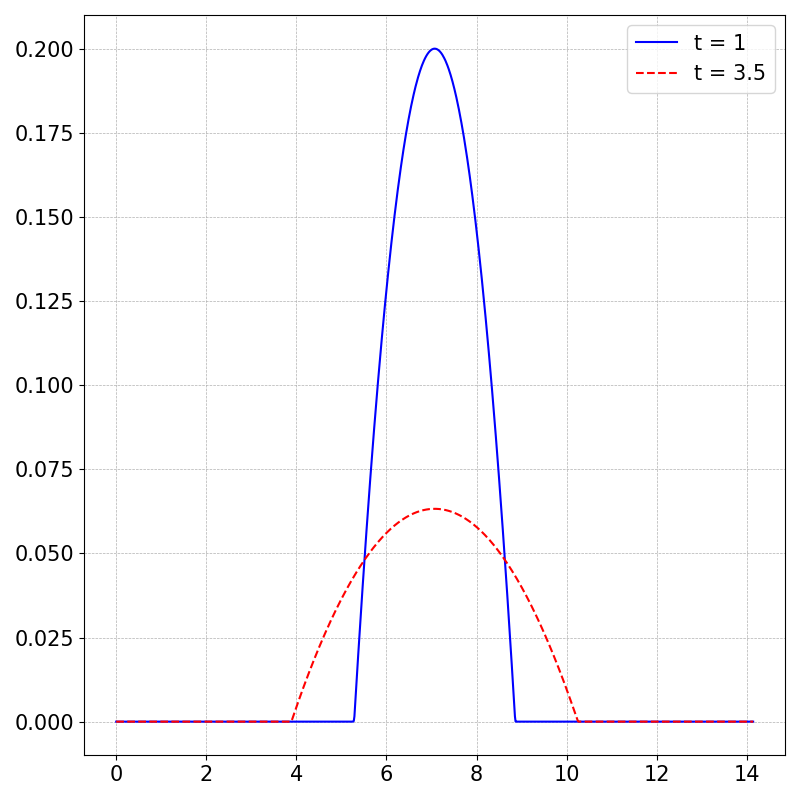}
         \caption{$\vec{v} = (0,0)^T$}
         \label{Upwind VS Wang V00}
     \end{subfigure}
     \hfill
     \begin{subfigure}[b]{0.4\textwidth}
         \centering
         \includegraphics[width=\textwidth]{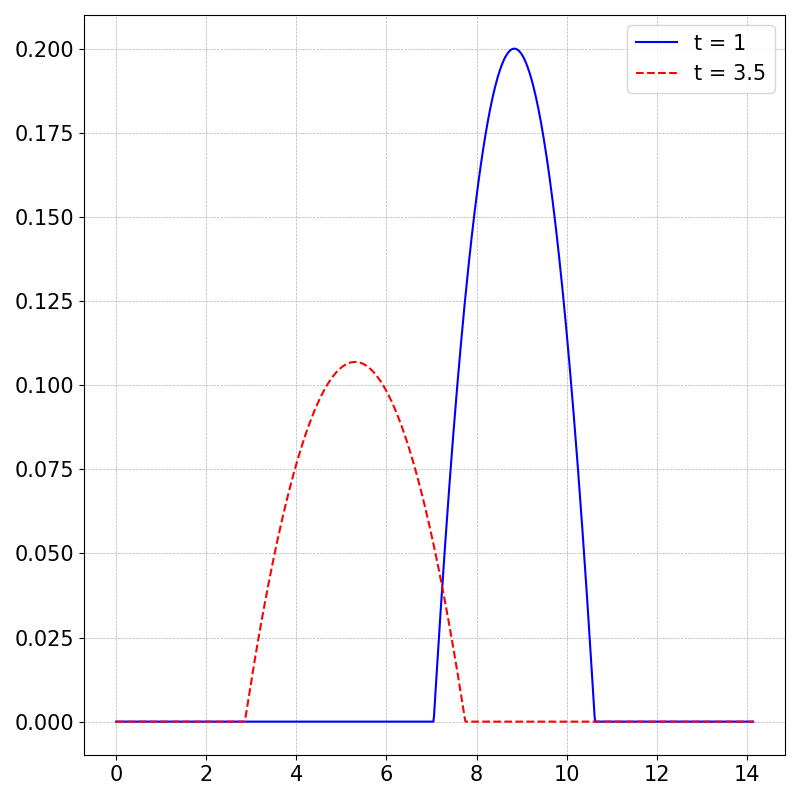}
         \caption{$\vec{v} = (1/2,1/2)^T$}
         \label{All schemes V00}
     \end{subfigure}
    \caption{The diagonal cross section of the analytical solution \eqref{analytical solution}.}
    \label{fig:corss section exact}
\end{figure}

Below we study the two cases $\vec{v} = (0, 0)^T$ and $\vec{v} = (1/2, 1/2)^T$ for which the analytical solutions are shown in Figures~\ref{fig:exact 2D} and \ref{fig:corss section exact}. The simulations are performed on the domain $\Omega = (-5, 5)^2$ over the time intervals $t \in (1, 10]$ for $\vec{v} = (0, 0)^T$ and $t \in (1, 7/2]$ for $\vec{v} = (0.5, 0.5)^T$.

For the temporal discretization, we combine the finite volume scheme with implicit Euler and the DG schemes with the second-order accurate diagonally Implicit Runge-Kutta (DIRK) method developed by Alexander \cite{alexander1977diagonally}. The spatial domain is discretized using an unstructured mesh with an initial mesh size of $h \approx 0.5$. The ratio of spatial mesh size and $\Delta t$ is kept constant.

The regularization coefficients for the cut-cell DG method are $\delta_1 = 2\cdot 10^{-5}$ and $\delta_2$ as detailed in Table~\ref{delta2_table}. Note that while the parameter $\delta_2$ could be further reduced, doing so does not decrease the overall numerical error. However, if $\delta_2$ is chosen too small, the nonlinear iteration may fail to converge. 
The initial value for the water height in dry regions is set to $\eta_0 = 4 \times 10^{-7}$, a value which does not affect the computational error.

{\renewcommand{\arraystretch}{1.1}
\setlength{\tabcolsep}{8pt}
\begin{table}
    \caption{The regularization coefficient ($\delta_2$) in each refinement level. \label{delta2_table}}
    \centering
    {\small\begin{tabular}{cccccc}
    \hline
        Refinement & 0 & 1 & 2 & 3 & 4 \\ \hline
        $\vec{v} = (0,0)^T$ & $1.0\cdot 10^{-3}$ & $1.0\cdot 10^{-3}$ & $1.0\cdot 10^{-3}$ & $1.0\cdot 10^{-3}$ & $3.5\cdot 10^{-4}$ \\ \hline
        $\vec{v} = (1/2,1/2)^T$ & $1.0\cdot 10^{-2}$ & $2.5\cdot 10^{-3}$ & $2.5\cdot 10^{-3}$ & $1.0\cdot 10^{-3}$ & $7.5\cdot 10^{-4}$ \\ \hline
    \end{tabular}}
\end{table}
}

\subsubsection{Discussion of convergence behavior}

{\renewcommand{\arraystretch}{1.05}
\setlength{\tabcolsep}{8pt}
\begin{table}
    \caption{$L^2$-error at final time and rate of convergence for the different schemes. Note that VFVM results are using $\mathbb{P}_1$ post-processing and implicit Euler time stepping while the DG schemes were combined with the second-order two-stage Alexander DIRK method.}
    \label{table:rate}
    \centering
     {\small\begin{tabular}{rrcccccc}
     \hline
     \hline
    $h$ & $\Delta t$ & DG-cutcell & RoC & DG & RoC & VFVM  & RoC \\ \hline
    \multicolumn{8}{c}{$\vec{v}=(0,0)^T$, $t=10$ } \\
    \hline
        1/2 & 1/2 & $9.41 \times 10^{-3}$ & ~ & $1.26 \times 10^{-2}$ & ~ & $3.39 \times 10^{-2}$ & ~ \\ 
        1/4 & 1/4 & $2.38 \times 10^{-3}$ & 1.98 & $3.53 \times 10^{-3}$ & 1.84 & $2.00 \times 10^{-2}$ & 0.76 \\ 
        1/8 & 1/8 & $5.84 \times 10^{-4}$ & 2.03 & $1.08 \times 10^{-3}$ & 1.71 & $1.22 \times 10^{-2}$ & 0.71 \\ 
        1/16 & 1/16 & $1.62 \times 10^{-4}$ & 1.85 & $3.75 \times 10^{-4}$ & 1.53 & $7.19 \times 10^{-3}$ & 0.77 \\ 
        1/32 & 1/32 & $4.86 \times 10^{-5}$ & 1.74 & $1.22 \times 10^{-4}$ & 1.62 & $4.07 \times 10^{-3}$ & 0.82 \\
      \hline
      \hline
    \multicolumn{8}{c}{$\vec{v}=(1/2,1/2)^T$, $t=7/2$}  \\
    \hline
        1/2 & 1/20 & $2.27 \times 10^{-2}$ & -- & $4.11 \times 10^{-2}$ & ~ & $1.55 \times 10^{-1}$ & ~ \\
        1/4 & 1/40 & $6.28 \times 10^{-3}$ & 1.86 & $1.68 \times 10^{-2}$ & 1.29 & $1.07 \times 10^{-1}$ & 0.54 \\
        1/8 & 1/80 & $1.68 \times 10^{-3}$ & 1.90 & $7.57 \times 10^{-3}$ & 1.15 & $6.87 \times 10^{-2}$ & 0.64 \\
        1/16 & 1/160 & $4.75 \times 10^{-4}$ & 1.82 & $3.59 \times 10^{-3}$ & 1.08 & $4.21 \times 10^{-2}$ & 0.71 \\
        1/32 & 1/320 & $1.09 \times 10^{-4}$ & 2.13 & $1.74 \times 10^{-3}$ & 1.04 & $2.49 \times 10^{-2}$ & 0.76 \\
     \hline
     \hline
    \end{tabular}}
\end{table}}

Table \ref{table:rate} shows the $L^2$-errors and convergence rates for the three different schemes and both cases. We observe that the VFVM with $P_1$-postprocessing reaches almost first order convergence for both cases asymptotically. The rates are worse on coarse meshes for the inclined plane. The convergence rate of the standard DG method is between $1.5$ and $2$ for the case $b=0$, while it is $1$ for the inclined plane. Finally, the cut-cell DG method attains order $2$ in both cases. Consequently, it reaches the smallest $L^2$-errors. On the inclined plain, the solution of the cut-cell DG method on the coarsest mesh matches the solution of the VFVM on a four-times refined mesh!    

\begin{figure}
     \centering
     \begin{subfigure}[b]{0.45\textwidth}
         \centering
         \includegraphics[width=\textwidth]{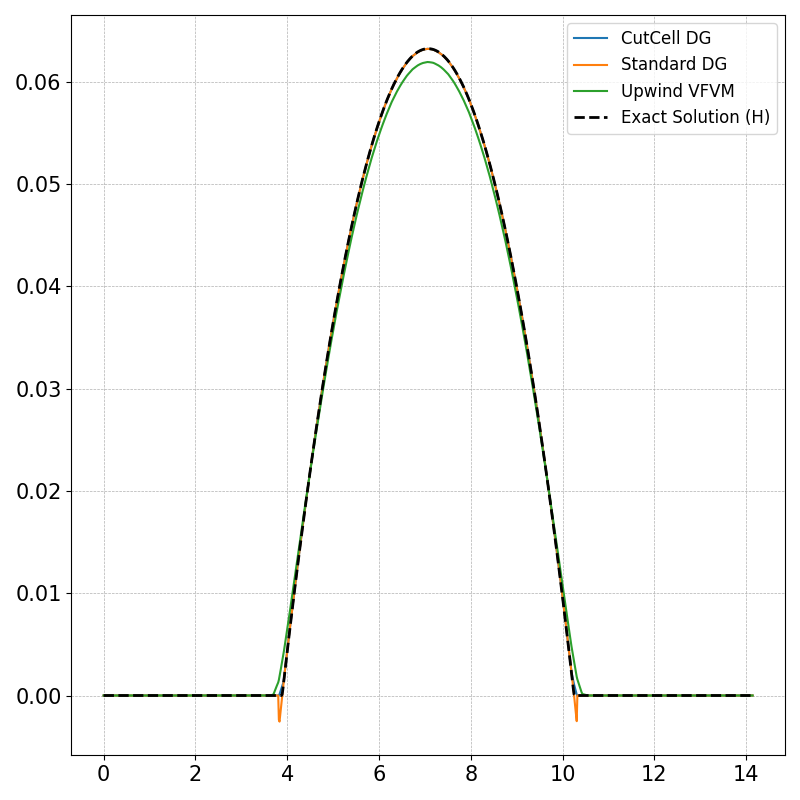}
         \caption{Overall domain}
         \label{All schemes V00}
     \end{subfigure}
     \hfill
     \begin{subfigure}[b]{0.45\textwidth}
         \centering
         \includegraphics[width=\textwidth]{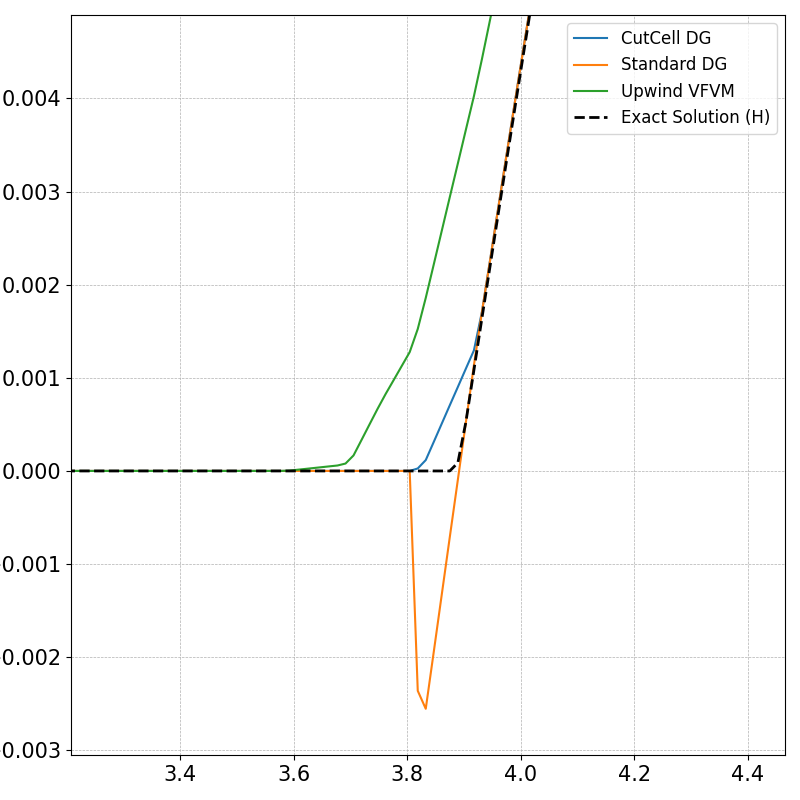}
         \caption{Wet/dry front at leading edge}
         \label{All schemes V00 ZoomIn}
     \end{subfigure}
    \caption{Diagonal cross section of the Barenblatt solution with $\vec{v} = 0$ using $h = 1/16$ and $\Delta t = 1/20$ using the three different schemes. The pulse moves from right to left.}
    \label{fig:corss section Barenblatt V00}
\end{figure}

\begin{figure}
     \centering
     \begin{subfigure}[b]{0.45\textwidth}
         \centering
         \includegraphics[width=\textwidth]{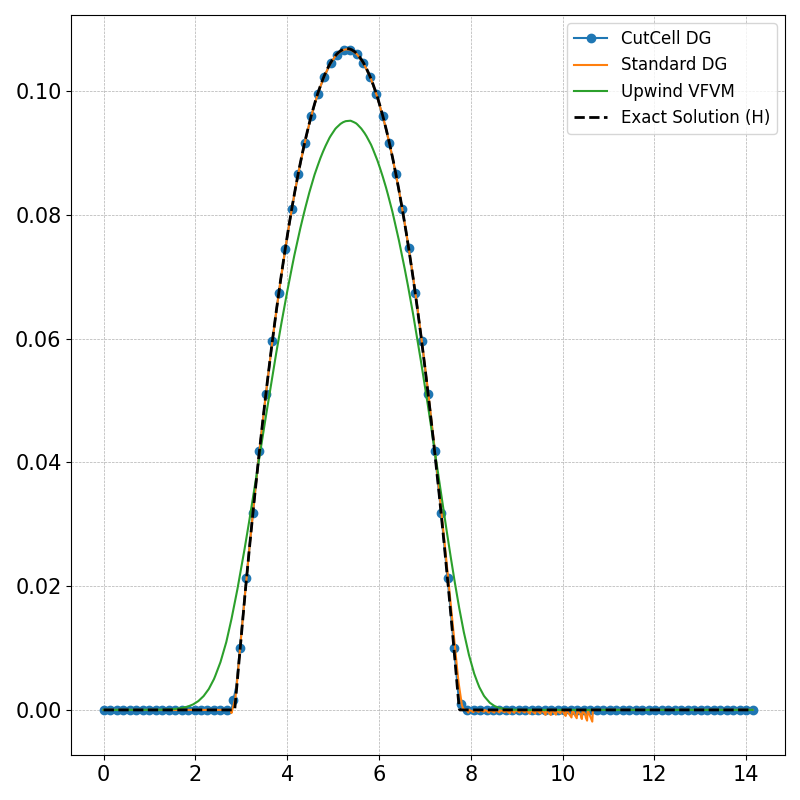}
         \caption{Overall Domain}
         \label{All scheme V55}
     \end{subfigure}
     \hfill
     \begin{subfigure}[b]{0.45\textwidth}
         \centering
         \includegraphics[width=\textwidth]{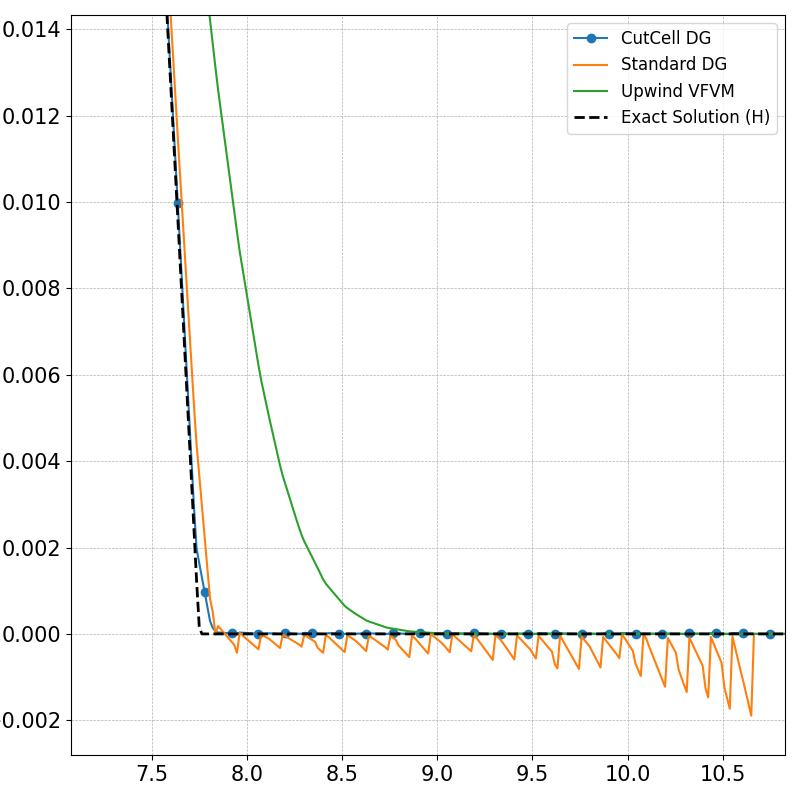}
         \caption{Wet/dry front at trailing edge}
         \label{All scheme V55 zoom}
     \end{subfigure}
    \caption{Diagonal cross section of the Barenblatt solution with $\vec{v} = (1/2,1/2)^T$, $h = 1/16$, $\Delta t = 1/320$ using the three different schemes. The pulse moves from right to left.}
    \label{V55}
\end{figure}


The reduced convergence rates of the VFVM and the standard DG scheme originate from the behaviour of the two schemes at the wet/dry-front. Figures \ref{fig:corss section Barenblatt V00} and \ref{V55} show cross sections along the diagonal of the domain with zooms of the leading wet/dry front (for $v=0$) and trailing wet/dry front (for $v=(1/2,1/2)^T$). 
Figure \ref{fig:corss section Barenblatt V00}(b) shows the leading front, meaning cells go from dry to wet state. VFVM has a much worse resolution of the wet/dry front compared to the DG schemes. The standard DG scheme has the best approximation of the wet/dry front at the cost of a negative water height. Given enough time, a cell is completely flooded and water height is positive in the whole cell. Cut-cell DG is very close to standard DG but never produces negative water height.
Figure \ref{V55}(b) shows the trailing front, i.e. cells go from wet to dry state. Again VFVM has a much worse approximation of the wet/dry front, with large error in position and shape. Standard DG has a good approximation of the position but many dry cells are left with negative water heights. This affects also the quality of the solution globally because this ``negative'' mass must be compensated by positive mass elsewhere. Clearly, the cut-cell DG scheme has an excellent approximation of the wet/dry front without negative water height.


\subsubsection{Discussion of computational cost}

\begin{figure}
    \centering
    \includegraphics[width=0.95\linewidth]{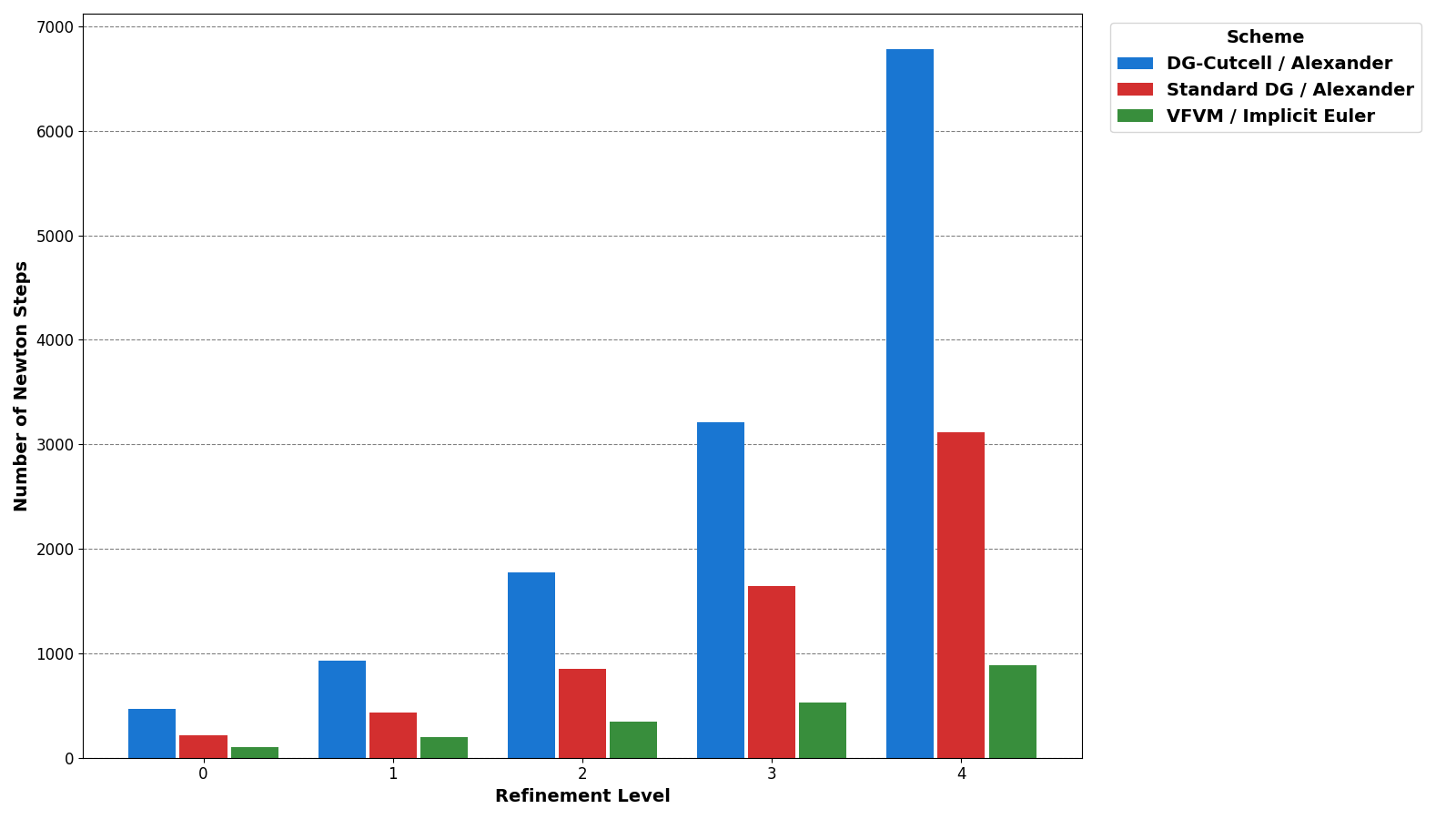}
    \caption{Total number of Newton iterations for the proposed schemes across five refinement levels (about 200000 triangular elements on level 4), using a velocity of $v = (1/2,1/2)^T$ and time step sizes as defined in the Table \ref{table:rate}.}
    \label{fig:newton_cost}
\end{figure}

\begin{figure}
    \centering
    \includegraphics[width=0.95\linewidth]{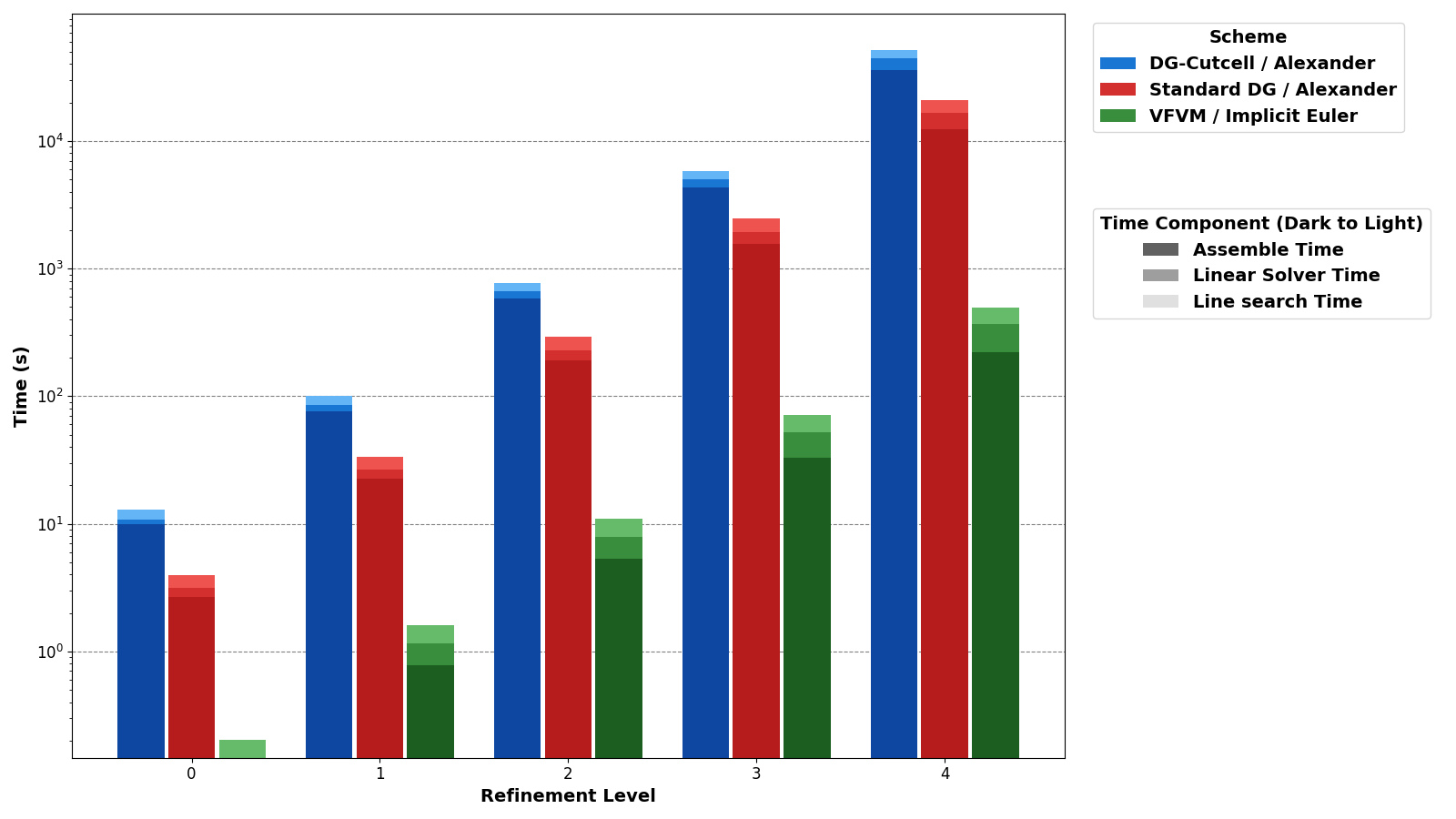}
    \caption{Total computation time for the proposed schemes across five refinement levels (about 200000 triangular elements on level 4), using a velocity of $v = (1/2,1/2)^T$ and time step sizes as defined in the Table \ref{table:rate}.}
    \label{fig:time_cost}
\end{figure}

The accuracy of a numerical scheme needs to be evaluated relative to its computational cost. Fig.~\ref{fig:newton_cost} shows the total number of Newton steps needed by the three schemes across five refinement levels, ranging from 800 triangular elements on level 0 to about 200000 elements on level 5 for the time step sizes given in Table \ref{table:rate}. Since the length of the time interval is fixed and time step size is reduced an increase of the number of Newton steps proportional to the number of time steps is expected. The number Newton steps is lowest for the VFVM and goes up by a factor three for the standard DG scheme. However, it has to be taken into account that VFVM is used with implicit Euler and the DG schemes are used with the two-stage Alexander scheme, so the number of nonlinear systems to be solved is doubled. Finally, the number of Newton steps is again doubled from the standard DG scheme to the cut-cell scheme.

Fig.~\ref{fig:time_cost} shows computational time where the schemes were run sequentially on one core of an AMD Epyc 7713 CPU. The increase in computational time over the levels reflects the scaling of temporal and spatial mesh size (factor of eight per refinement). The standard DG scheme takes about 30 times longer than the VFVM and the cut-cell DG scheme is about 2.5 times more expensive than the standard DG scheme for the same discretization parameters. Note, however, that our implementation of the cut-cell DG scheme was not optimized for performance yet (nor was the finite volume scheme).

Although the DG schemes are considerably more computationally expensive for the same discretization parameters, taking into account the accuracy of the solution changes the picture. According to Table \ref{table:rate} the cut-cell DG scheme on level 1 is more accurate in both cases than the VFVM solution on level 4 and therefore cheaper to compute.


\subsection{Obstacle example } \label{sec:obstacle}

\begin{figure}
     \centering
     \begin{subfigure}[b]{0.5\textwidth}
         \centering
         \includegraphics[width=\textwidth]{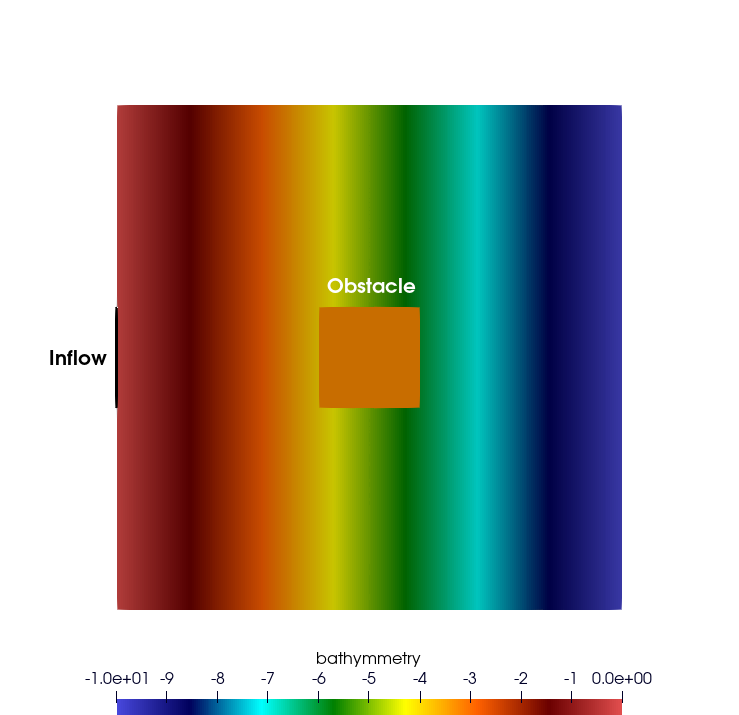}
         \caption{Bathymmetry}
         \label{domain setup obstacle}
     \end{subfigure}
     \hfill
     \begin{subfigure}[b]{0.45\textwidth}
         \centering
         \includegraphics[width=\textwidth]{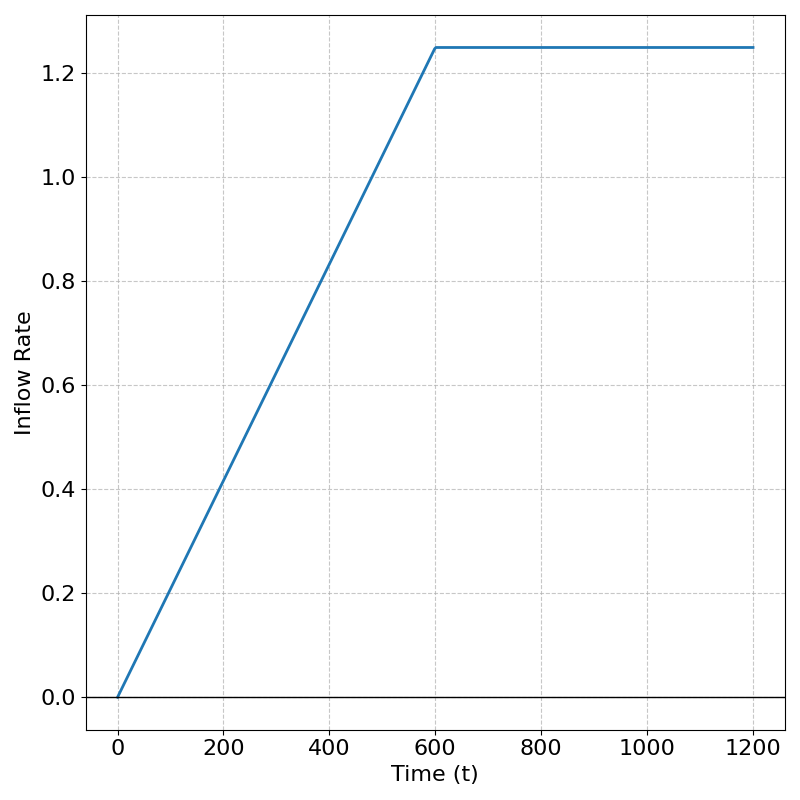}
         \caption{Inflow rate}
         \label{Inflow rate}
     \end{subfigure}
    \caption{The setup domain includes bathymetry, obstacle position, inflow boundary, and inflow rate.}
    \label{obstacle example setting}
\end{figure}

In this section, we use Manning's friction law with $n = 0.04$, $\alpha=5/3$ and $\gamma=1/2$:
\begin{equation}
    \partial_t u(x,t) - 25 \nabla \cdot \left(\frac{H(x,t)^{5/3}(x,t)}{G(x,t)^{1/2}}\nabla u \right) = 0.
\end{equation}
The simulation is performed on the domain $\Omega = (0,1000)^2$ using an unstructured mesh with an element size of $h \approx 10~\text{m}$. The obstacle is located at the box $\mathcal{B} = \{ (x_0,x_1)^T \,:\, 400 < x_0,x_1 < 600\}$. The bathymmetry is given by
\begin{equation}
    b(x) = \begin{cases}
        -0.01x_0, &x \in \Omega \setminus \mathcal{B}, \\
        -3.25, &x \in \mathcal{B}.
    \end{cases}
\end{equation}
A homogeneous Neumann boundary condition is applied to all boundaries except at the inflow boundary, which is illustrated in Fig.~\ref{domain setup obstacle}. The inflow rate is specified as shown in Fig.~\ref{Inflow rate}, time step size is $\Delta t = 10 \, s$ and final time is $1200~\text{s}$. Here, we apply the Alexander DIRK method in time for both CUDGM and VFVM. The regularization parameters in CUDGM are set as $\delta_1 = 2 \cdot 10^{-5}$, $\delta_2 = 1 \cdot 10^{-3}$, and $\eta_0 = 1 \cdot 10^{-5}$. The initial water height of VFVM is set to be zero.

\begin{figure}
     \centering
     \begin{subfigure}[b]{0.45\textwidth}
         \centering
         \includegraphics[width=\textwidth]{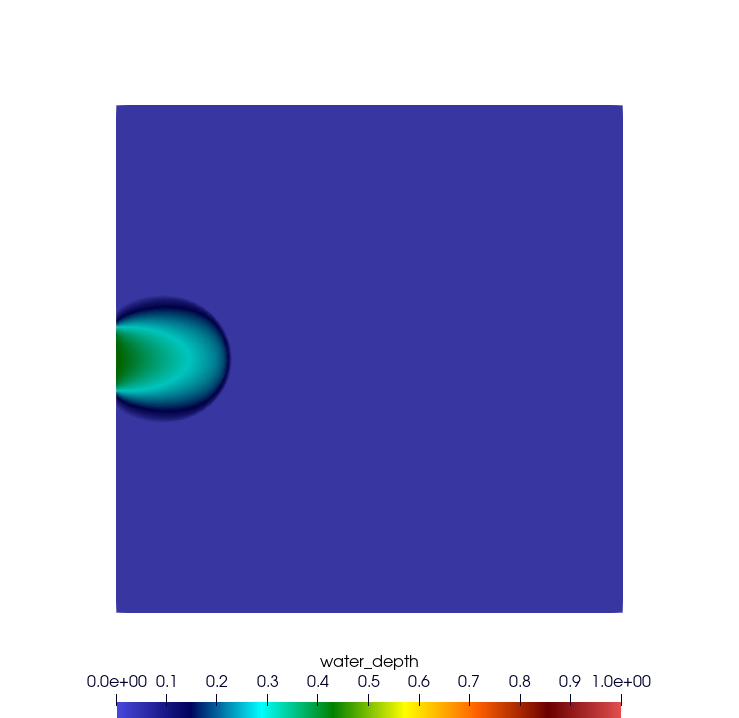}
         \caption{$t = 300$}
         \label{obstacle_cutcell__refine1_t300}
     \end{subfigure}
     \hfill
     \begin{subfigure}[b]{0.45\textwidth}
         \centering
         \includegraphics[width=\textwidth]{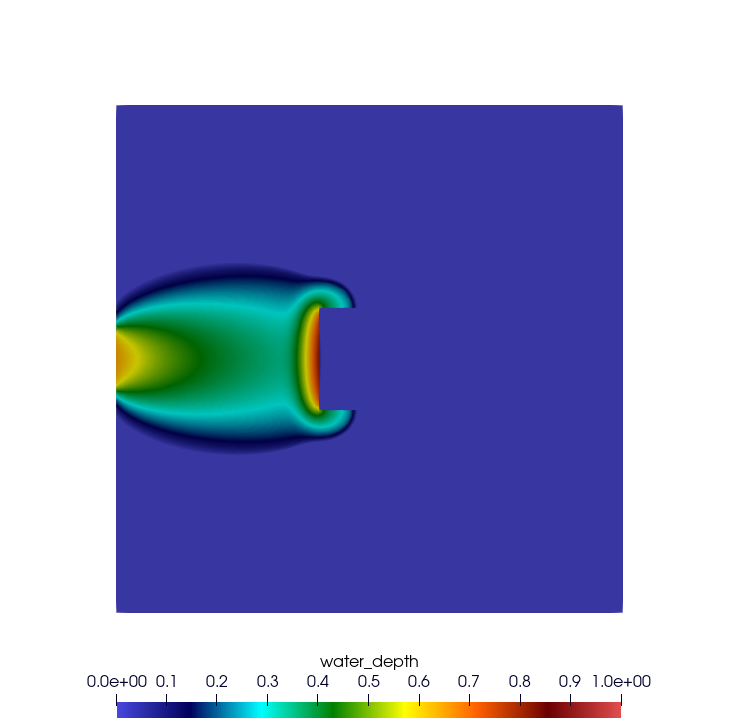}
         \caption{$t = 600$}
         \label{obstacle_cutcell__refine1_t600}
     \end{subfigure}
     \begin{subfigure}[b]{0.45\textwidth}
         \centering
         \includegraphics[width=\textwidth]{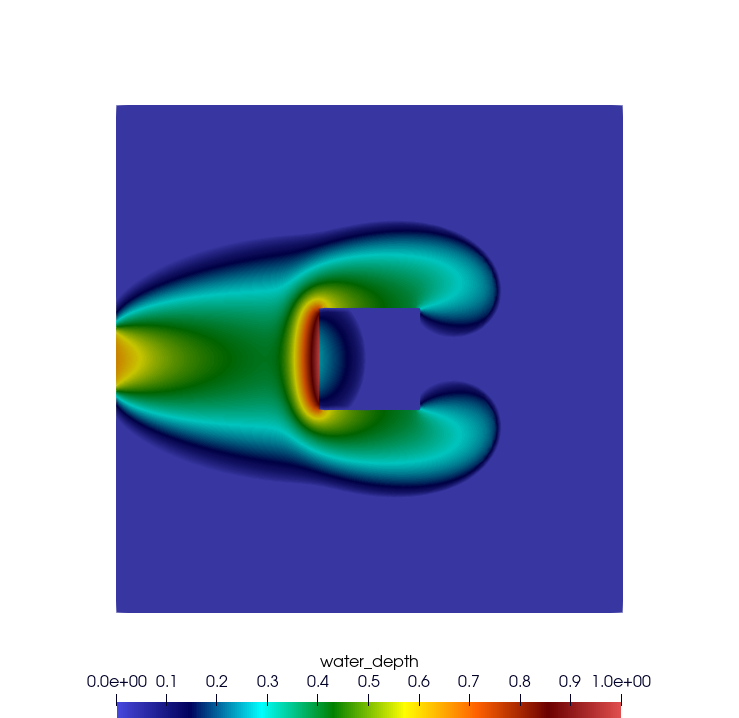}
         \caption{$t = 900$}
         \label{Cutcell_DG_obstacle_StructuredMesh_FullModel}
     \end{subfigure}
     \hfill
     \begin{subfigure}[b]{0.45\textwidth}
         \centering
         \includegraphics[width=\textwidth]{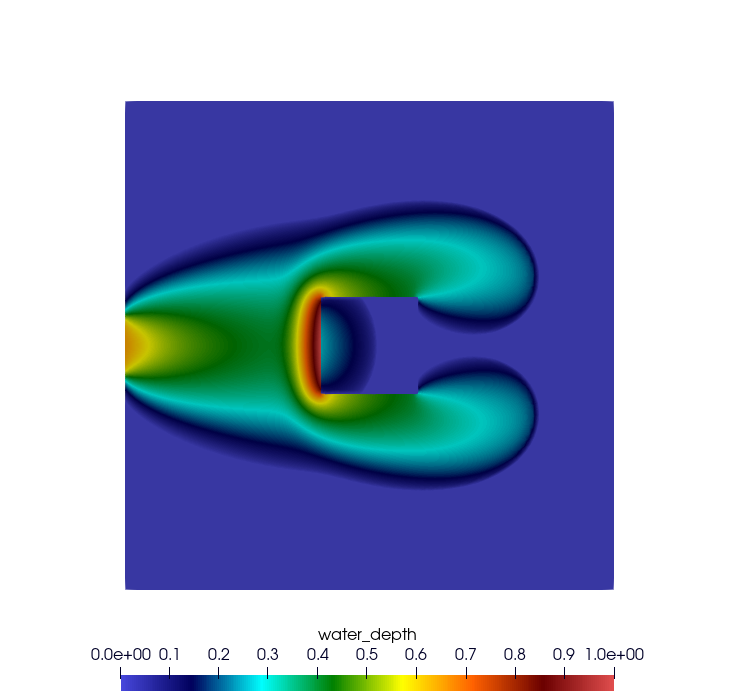}
         \caption{$t = 1000$}
         \label{Cutcell_DG_obstacle_UnstructuredMesh_FullModel}
     \end{subfigure}
    \caption{The numerical solution of the flow on the inclined plain with obstacle implemented by cut-cell DG method using $h \approx 10$ and $\Delta t = 10$.}
    \label{obstacle example 2D result}
\end{figure}

\begin{figure}
    \centering
    \includegraphics[width=1.0\textwidth]{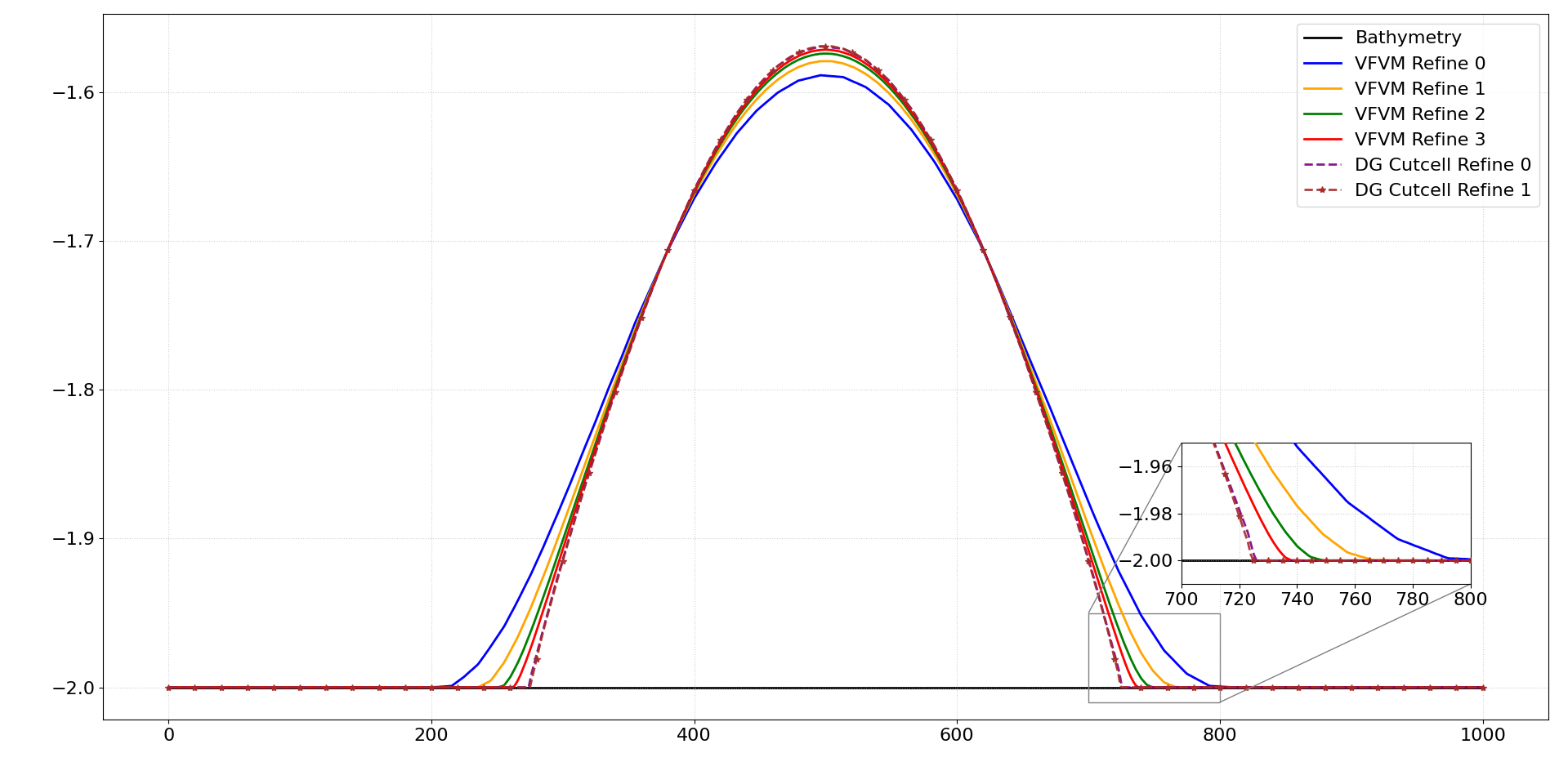}
    \caption{Cross-section of the solutions on the axis $x_0 = 200$ at $t = 1000$.}
    \label{fig:crossection_obstacle_x}
\end{figure}

\begin{figure}
    \centering
    \includegraphics[width=1.0\textwidth]{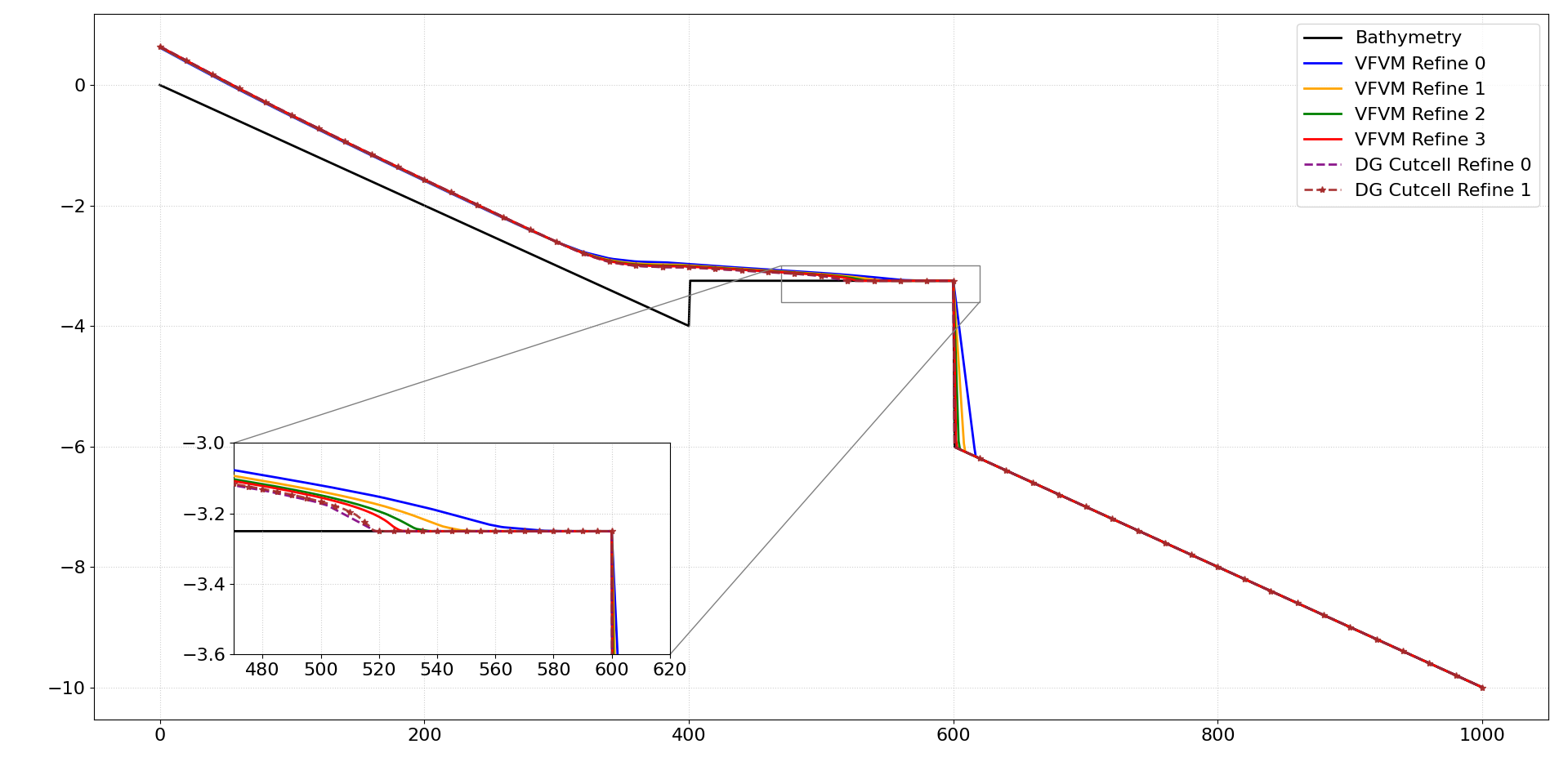}
    \caption{Cross-section of the solutions on the axis $x_1 = 500$ at $t = 1000$.}
    \label{fig:crossection_obstacle_y}
\end{figure}

Figure~\ref{obstacle example 2D result} illustrates the time evolution of the 2D flow simulation. A water wave propagates from the inflow boundary, impinging upon the obstacle at $t \approx 600~\mathrm{s}$. The flow is diverted around the structure, and by $t = 900~\mathrm{s}$, significant inundation is observed within the obstacle's interior. Subsequently, the main wavefront continues to propagate towards the right boundary.

To observe the convergence, we analyzed cross-sections at $x_0=200$ and $x_1=500$ at $t = 1000~\mathrm{s}$, shown in Fig.~\ref{fig:crossection_obstacle_x} and Fig.~\ref{fig:crossection_obstacle_y}, respectively. These plots include a zoom-in window at the wet/dry front to illustrate the interface behavior. The results show that the VFVM solution converges to the CUDGM solution, and the CUDGM results at different refinement levels are not noticeably different. The CUDGM solution on the coarsest mesh level is more accurate than the VFVM solution after three refinements. Also note that the steep slope of the water level in the VFVM in Fig.~\ref{fig:crossection_obstacle_y} at $x_0 = 600$ is actually due to a piecewise linear graphical representation of the bathymmetry in VFVM.

\subsection{Dambreak in the river valley} \label{Sec:dam_break}

This example also employs the full DWE model with Manning's friction law and Manning's coefficient $n = 0.04$. The example serves to demonstrate that CUDGM can be applied to a realistic problem with relatively complicated bathymmetry. Furthermore, it provides a practical illustration of the non-negativity theorem of the VFVM.

\subsubsection{Data description}

\begin{figure}[ht]
    \centering
    \hspace*{-0.5in}
    \includegraphics[scale=0.3]{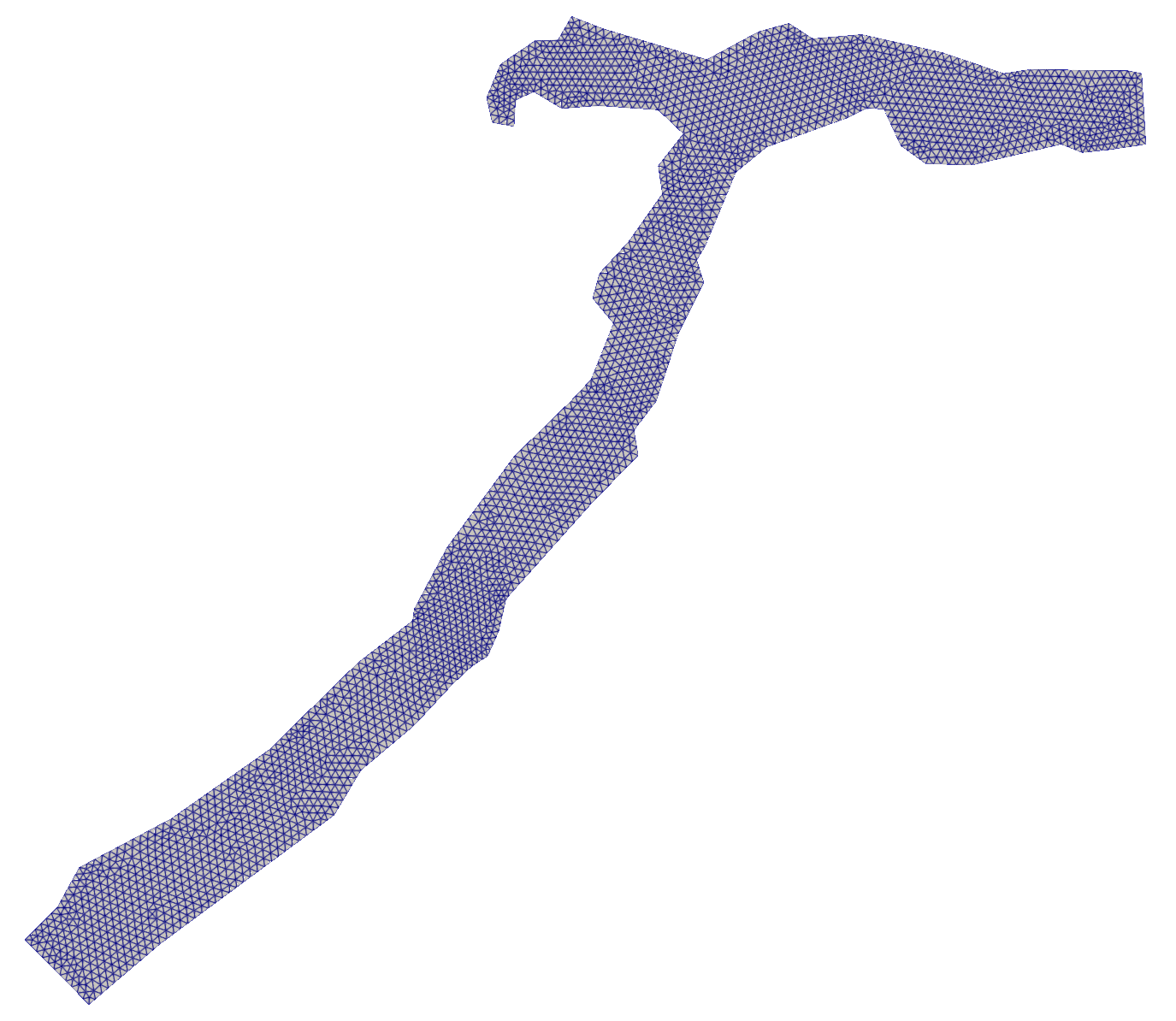}
    \caption{The unstructured triangular mesh applied to the domain from \cite{neelz2010benchmarking}}
    \label{fig:downflow mesh}
\end{figure}

\begin{figure}[ht]
    \centering
    \hspace*{-0.5in}
    \includegraphics[scale=0.3]{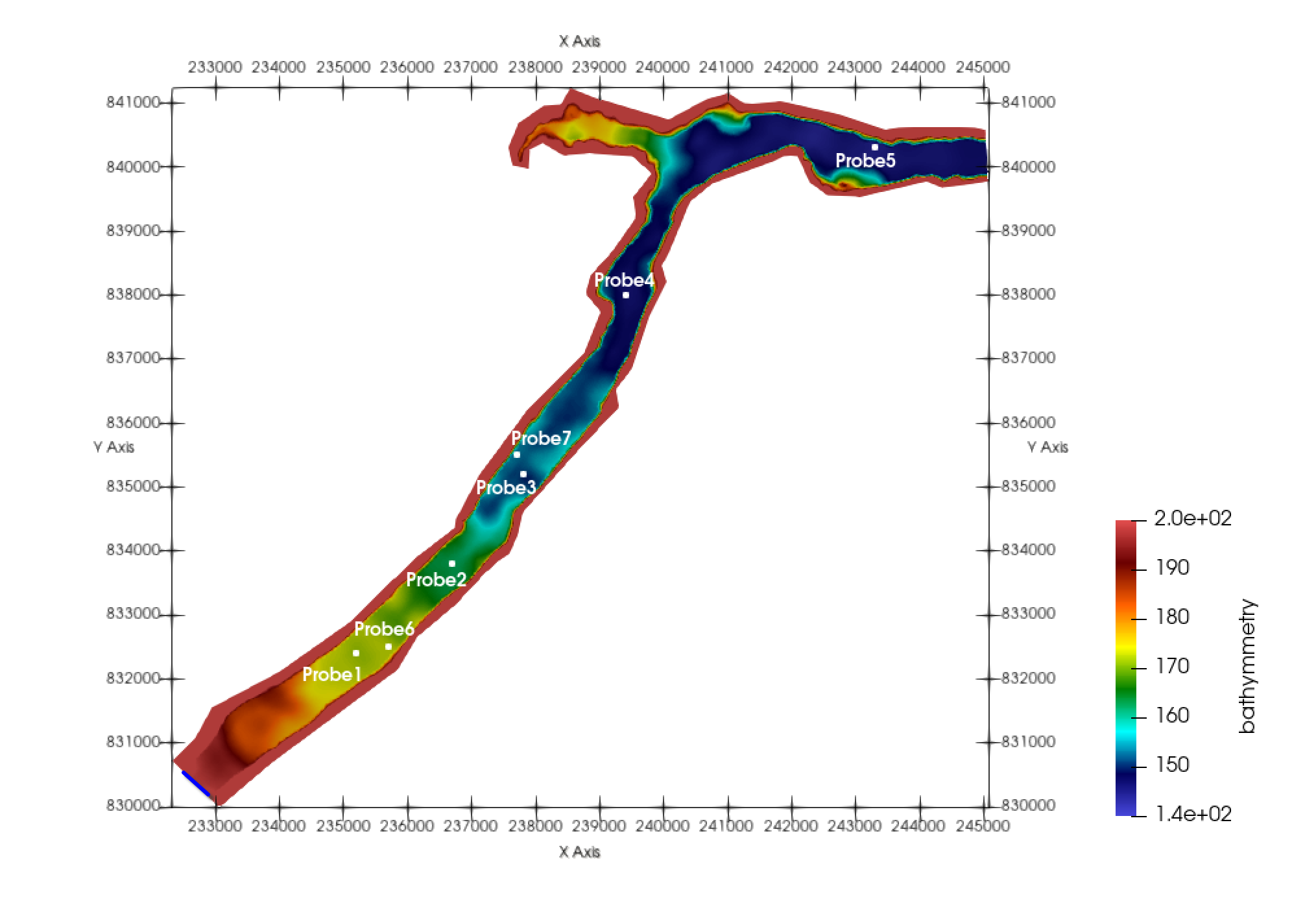}
    \caption{Bottom elevation and stations for measuring water depth from \cite{neelz2010benchmarking}}
    \label{fig:bathymmetry_test5}
\end{figure}

\begin{figure}[ht]
    \centering
    \hspace*{-0.5in}
    \includegraphics[scale=0.4]{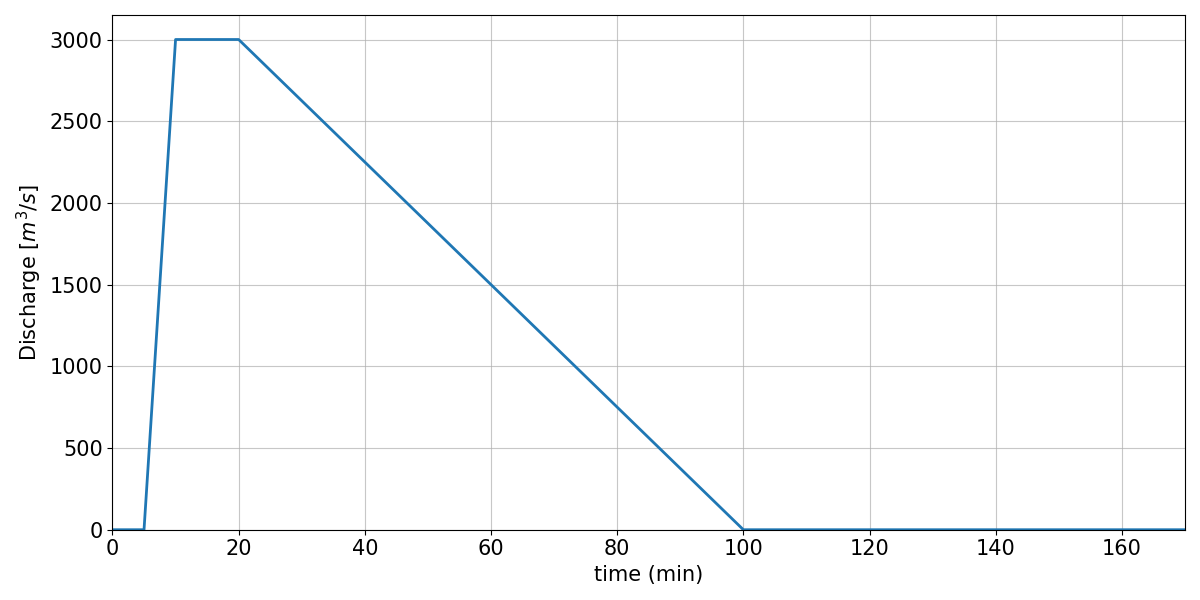}
    \caption{The inflow rate from \cite{neelz2010benchmarking}}
    \label{fig:inflow_rate}
\end{figure}

We use a model benchmark from \cite{neelz2010benchmarking}, which provides the domain geometry, bathymetry, and inflow conditions. The model covers a $19.1\,\text{km}^2$ area, with bathymetry defined by a $10\,\text{m}$ resolution raster dataset. The domain is discretized using an unstructured triangular mesh with an element size of $h \approx 132\,\text{m}$, resulting in $5442$ total elements. The coarse unstructured mesh is provided in Fig.~\ref{fig:downflow mesh}.

Fig~\ref{fig:bathymmetry_test5} shows the domain's surface elevation and the locations of measurement probes. A zero Neumann condition is enforced on all boundaries except for the specified inflow boundary, which is the blue line located on the lower left of the domain. The inflow rate is plotted in Fig.~\ref{fig:inflow_rate}.

The bathymetry data from the $10\,\text{m}$ resolution raster dataset is incorporated into the unstructured mesh by assigning the value of the nearest raster grid point to each mesh node. Subsequently, this nodal data is used for interpolation within the elements: a piecewise linear continuous interpolation is applied for the CUDGM, and a piecewise constant interpolation on Voronoi cells is used for the VFVM. To establish the bathymetry for the subsequent, finer-resolution meshes, we apply linear interpolation based on the initial coarse grid. This approach ensures a consistent and fixed bathymetry dataset for robust model comparison across different mesh resolutions.

\subsubsection{Simulation setting}
The simulation models a dam break scenario over a total duration of $54000\,\text{s}$ (15 hours). The inflow rate, shown in Figure~\ref{fig:inflow_rate}, rises rapidly at $t = 5\,\text{min}$, maintains a constant high rate, and then gradually decreases until $t = 100\,\text{min}$. An adaptive time-stepping scheme is employed, initialized at $\Delta t_0 = 1.0\,\text{s}$ and capped at $\Delta t_{\max} = 20.0\,\text{s}$. The step size is dynamically adjusted by a factor of $\sqrt{2}$ based on the Newton solver's performance: the time step is reduced upon non-convergence, or increased after three consecutive successful Newton iterations.

In the CUDGM, the regularization parameters vary by refinement level. At level 0, we set $\delta_1 = 2 \cdot 10^{-3}$ and $\delta_2 = 0.25$, with an initial condition of $\eta_0 = 1 \cdot 10^{-3}$. Subsequently, at refinement level 1, these values are adjusted to $\delta_1 = 2 \cdot 10^{-5}$, $\delta_2 = 0.1$, and $\eta_0 = 1 \cdot 10^{-5}$. It is important to note that while the CUDGM necessitates this positive initial water depth, the VFVM can initialize with zero water height, effectively simulating a dry surface.

\subsubsection{Result discussion}

\begin{figure}[h!]
     \centering
     \begin{subfigure}[b]{0.45\textwidth}
         \centering
         \includegraphics[width=\textwidth]{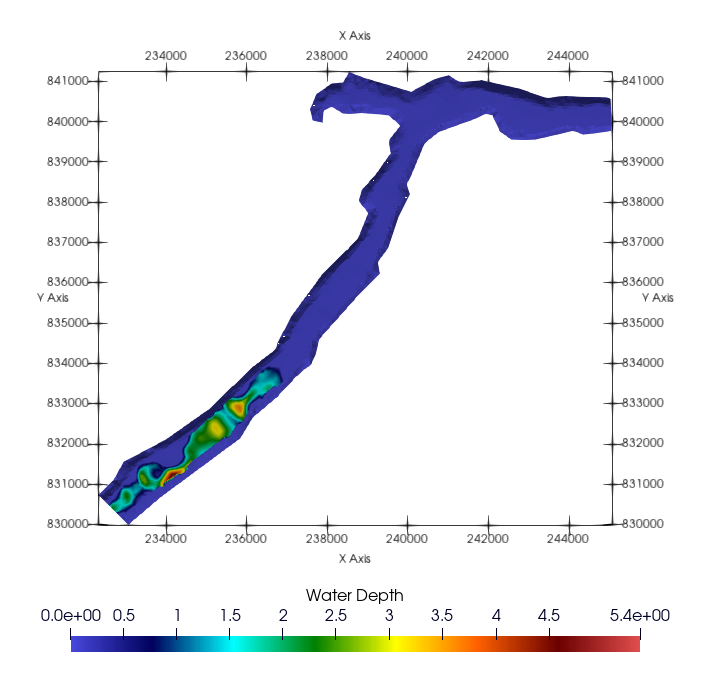}
         \caption{$t = 2500 $}
         \label{2d_result_t2500}
     \end{subfigure}
     \hfill
     \begin{subfigure}[b]{0.45\textwidth}
         \centering
         \includegraphics[width=\textwidth]{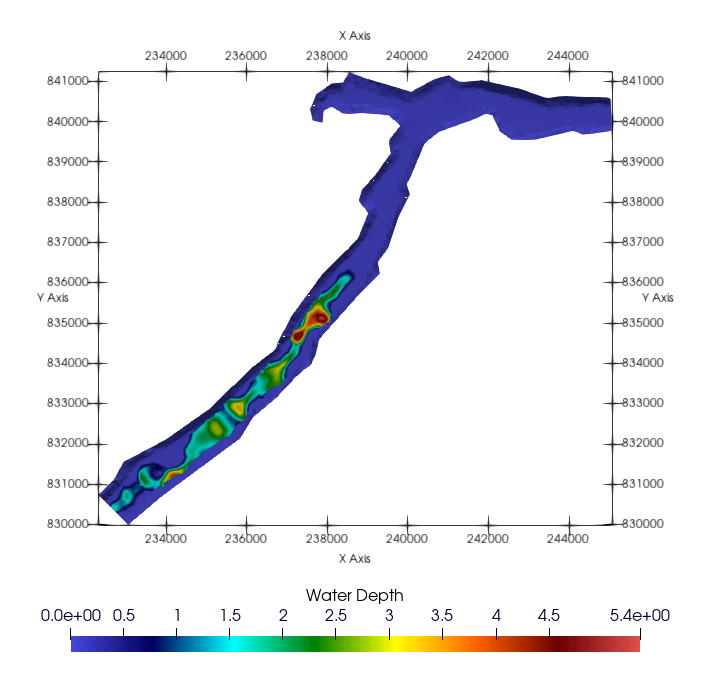}
         \caption{$t = 4000 $}
         \label{2d_result_t4000}
     \end{subfigure}
     \begin{subfigure}[b]{0.45\textwidth}
         \centering
         \includegraphics[width=\textwidth]{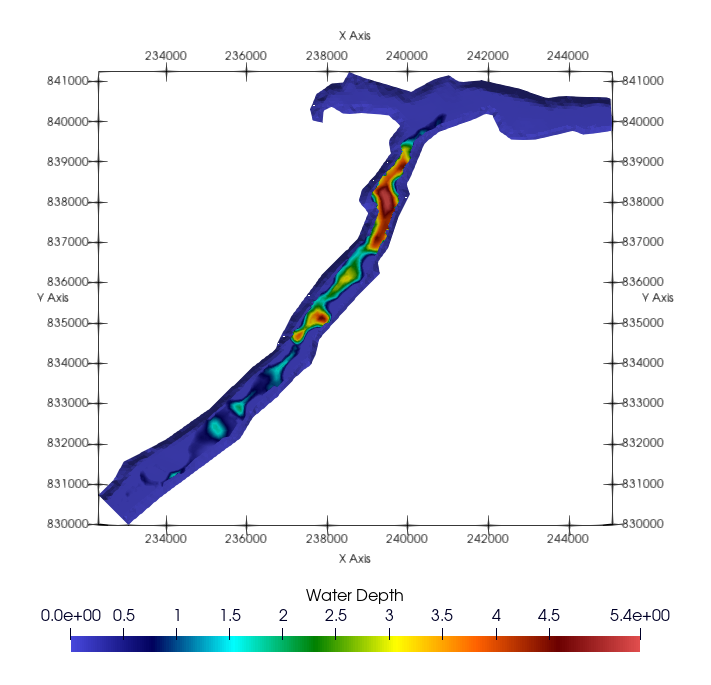}
         \caption{$t = 8000$}
         \label{2d_result_t8000}
     \end{subfigure}
     \hfill
     \begin{subfigure}[b]{0.45\textwidth}
         \centering
         \includegraphics[width=\textwidth]{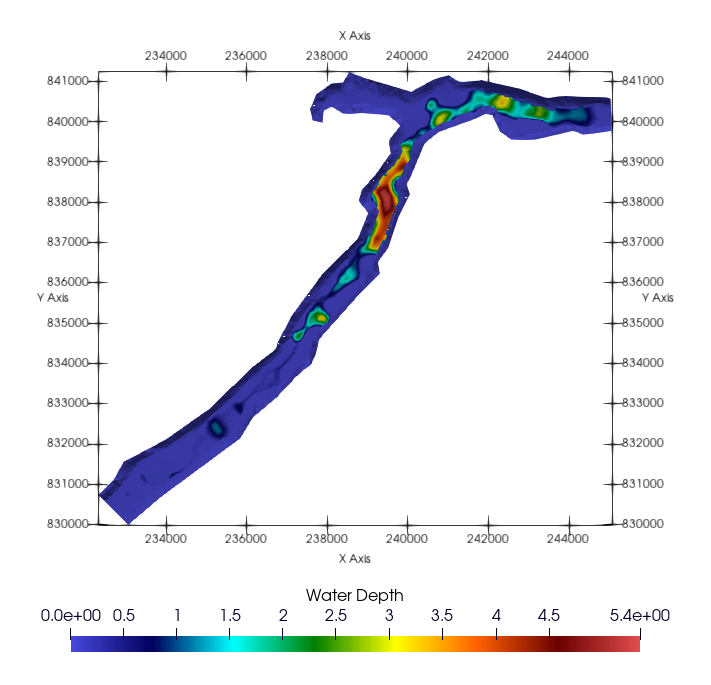}
         \caption{$t=16000$}
         \label{2d_result_t16000}
     \end{subfigure}
    \caption{The simulation of the water flow in the river valley using an unstructured mesh using the proposed cut cell-DG scheme with $h \approx 66 m$.}
    \label{down_folw_2d_simulation}
\end{figure}

Fig~\ref{down_folw_2d_simulation} illustrates the top-down view of the 2D simulation at several time steps. Initially, the flow velocity is high, as indicated by the significant evolution of the flow field between $t=2500~\text{s}$ and $t=4000~\text{s}$. After the inflow is halted at $t=6000~\text{s}$, the flow decelerates substantially. 
The drying process commences from the inflow region around $t=8000~\text{s}$, and by $t=16000~\text{s}$, the upstream half of the domain is almost completely dry.

\begin{table}[tbp]
\centering
\begin{tabular}{lcccc}
\toprule
Refinement & 0 & 1 & 2 & 3 \\
\midrule
$\min H$ & $7.59 \cdot 10^{-9}$ & $6.00 \cdot 10^{-9}$ & $4.10 \cdot 10^{-9}$ & $2.35 \cdot 10^{-9}$ \\
\bottomrule
\end{tabular}
\caption{Minimum water depth of the VFVM solution.}
\label{tab:minH_VFVM} 
\end{table}

To confirm the non-negativity preserving property of VFVM, as provided in Theorem \ref{thm:nonnegativity_VFVM}, we compute the minimum water depth over the complete simulation by $\min H = \min_{x \in \mathcal{X}_h} \min_{t \in \{t_i\}_{i = 0}^N} H(x,t)$. The results are shown in Table \ref{tab:minH_VFVM}. The minimum water depth stays positive, even though it is smaller in higher refinement levels. Note that CUDGM is positive preserving by definition of the solution space.

\begin{figure}
    \centering
    \includegraphics[width=0.8\linewidth]{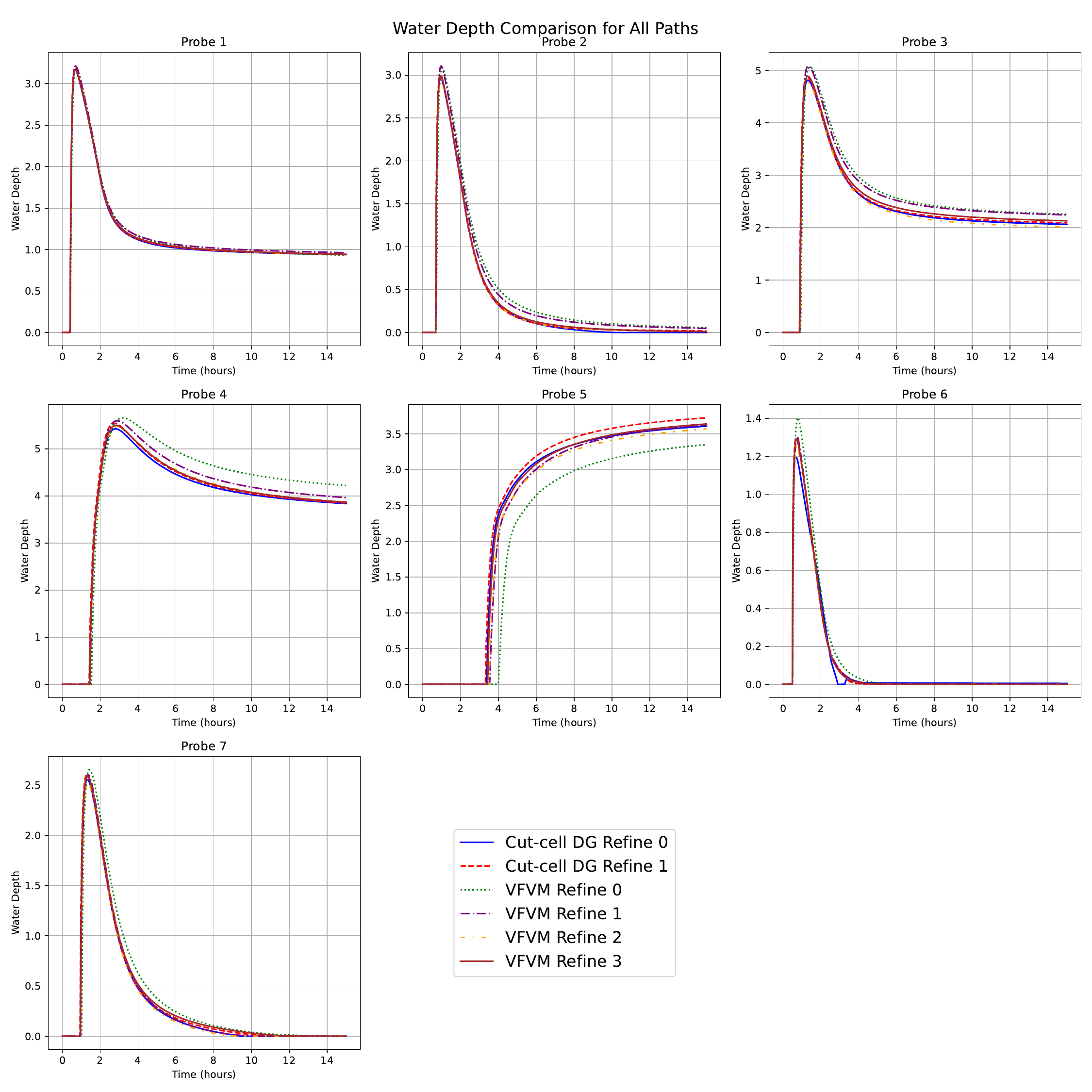}
    \caption{The hydrograph of the water depth in each probe.}
    \label{fig:water_depth_test5}
\end{figure}

We now present hydrographs (water depth over time) at seven different locations given in Table \ref{tab:probe_pos} and Figure \ref{fig:bathymmetry_test5}. Water depth is measured differently for each method: For the cut-cell DG solution, the depth is evaluated exactly at the probe locations which are all within a triangle. For the VFVM solution, the measurement is taken at the nearest computational point to the probe, which is equivalent to measuring the water depth in the $W^0_h$ space.

\begin{table}[tbp]
\centering
\begin{tabular}{lccccccc}
\toprule
Point ID & 1 & 2 & 3 & 4 & 5 & 6 & 7 \\
\midrule
$x_0$ & 235200 & 236700 & 237800 & 239400 & 243300 & 235700 & 237700 \\
$x_1$ & 832400 & 833800 & 835200 & 838000 & 840300 & 832500 & 835500 \\
\bottomrule
\end{tabular}
\caption{The position of the probes. \label{tab:probe_pos}}
\end{table}

\begin{figure}
     \centering
     \begin{subfigure}[b]{0.45\textwidth}
         \centering
         \includegraphics[width=\textwidth]{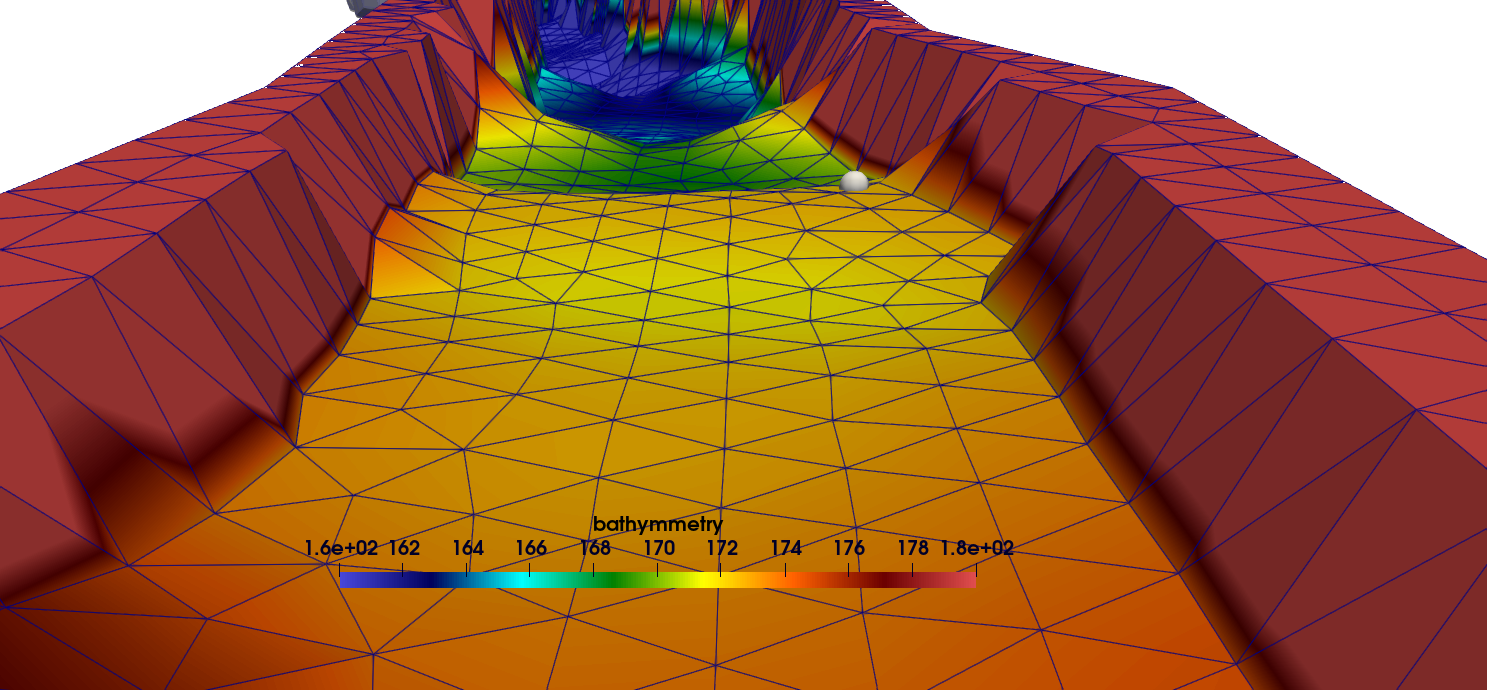}
         \caption{$t = 0 [s]$}
         \label{3d_result_t0}
     \end{subfigure}
     \hfill
     \begin{subfigure}[b]{0.45\textwidth}
         \centering
         \includegraphics[width=\textwidth]{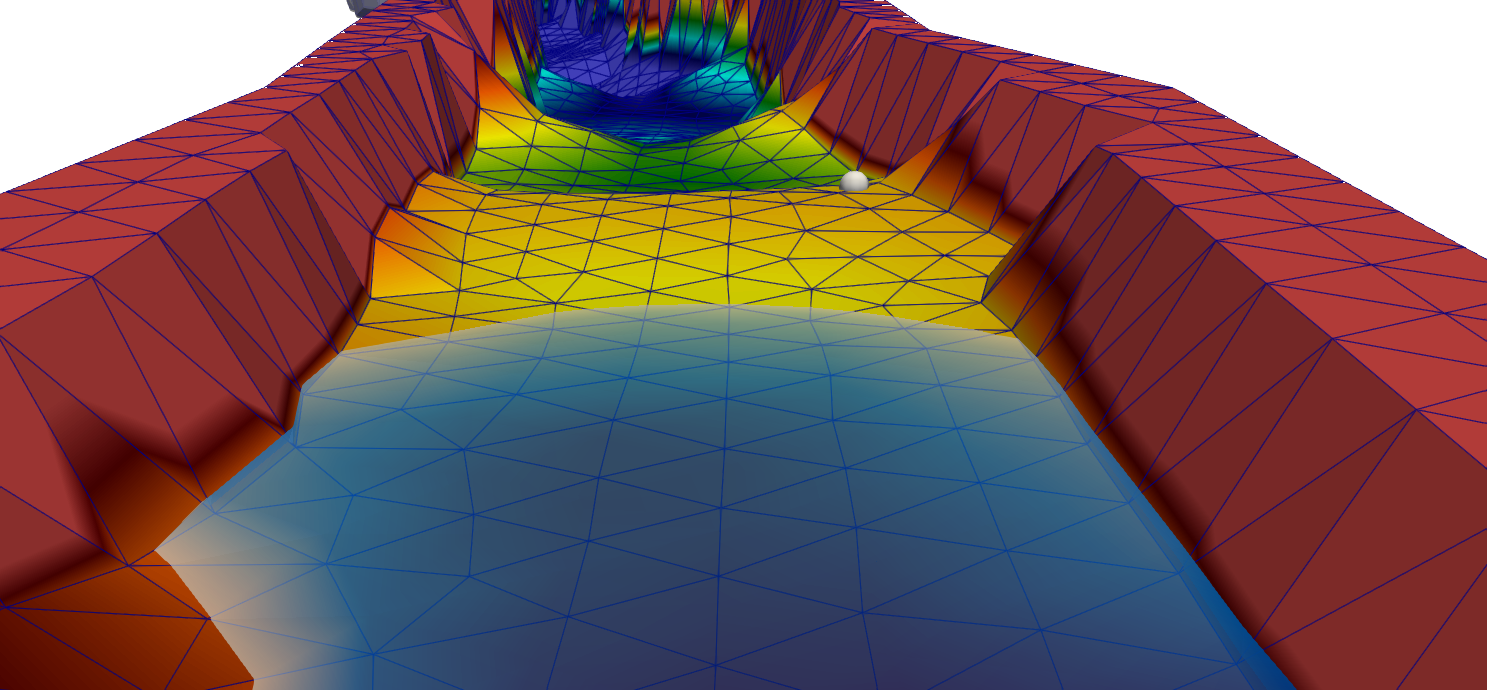}
         \caption{$t = 1500 $}
         \label{3d_result_t1500}
     \end{subfigure}
     \begin{subfigure}[b]{0.45\textwidth}
         \centering
         \includegraphics[width=\textwidth]{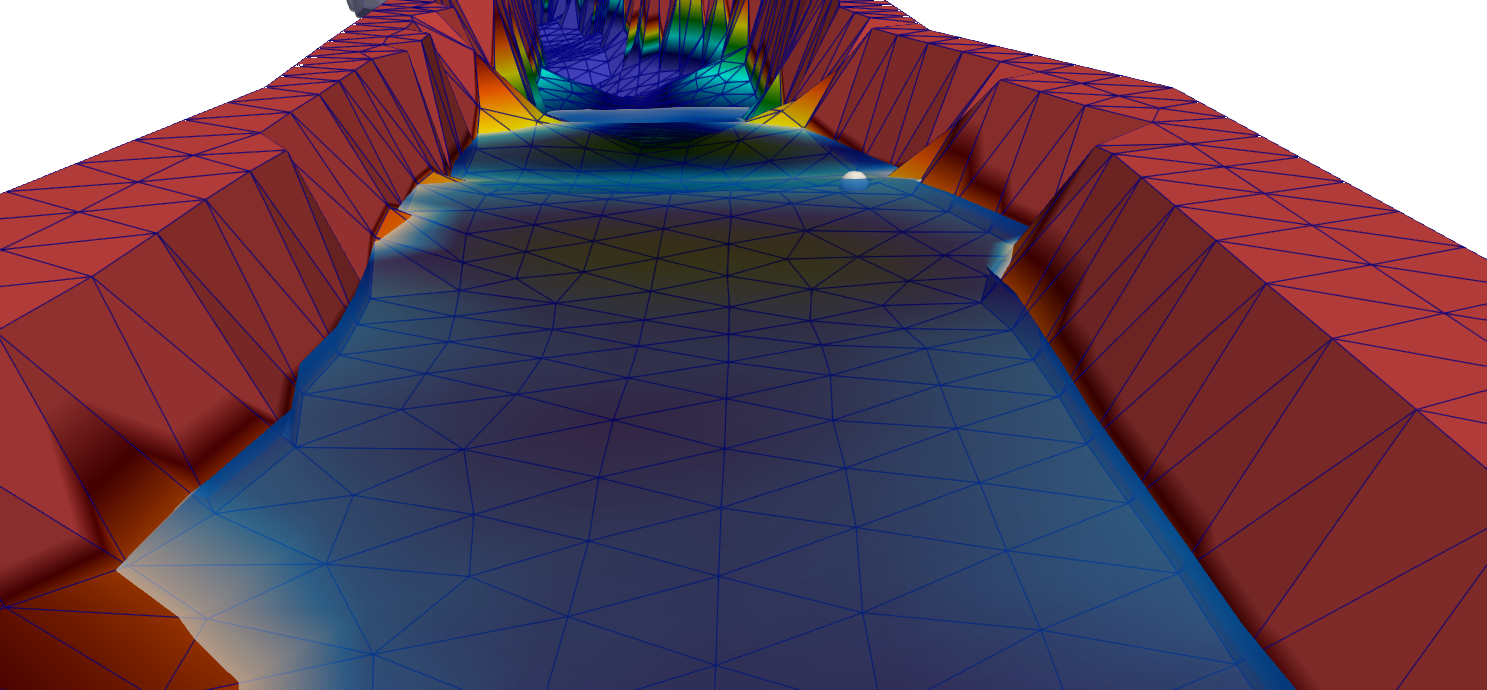}
         \caption{$t = 2500$}
         \label{3d_result_t2500}
     \end{subfigure}
     \hfill
     \begin{subfigure}[b]{0.45\textwidth}
         \centering
         \includegraphics[width=\textwidth]{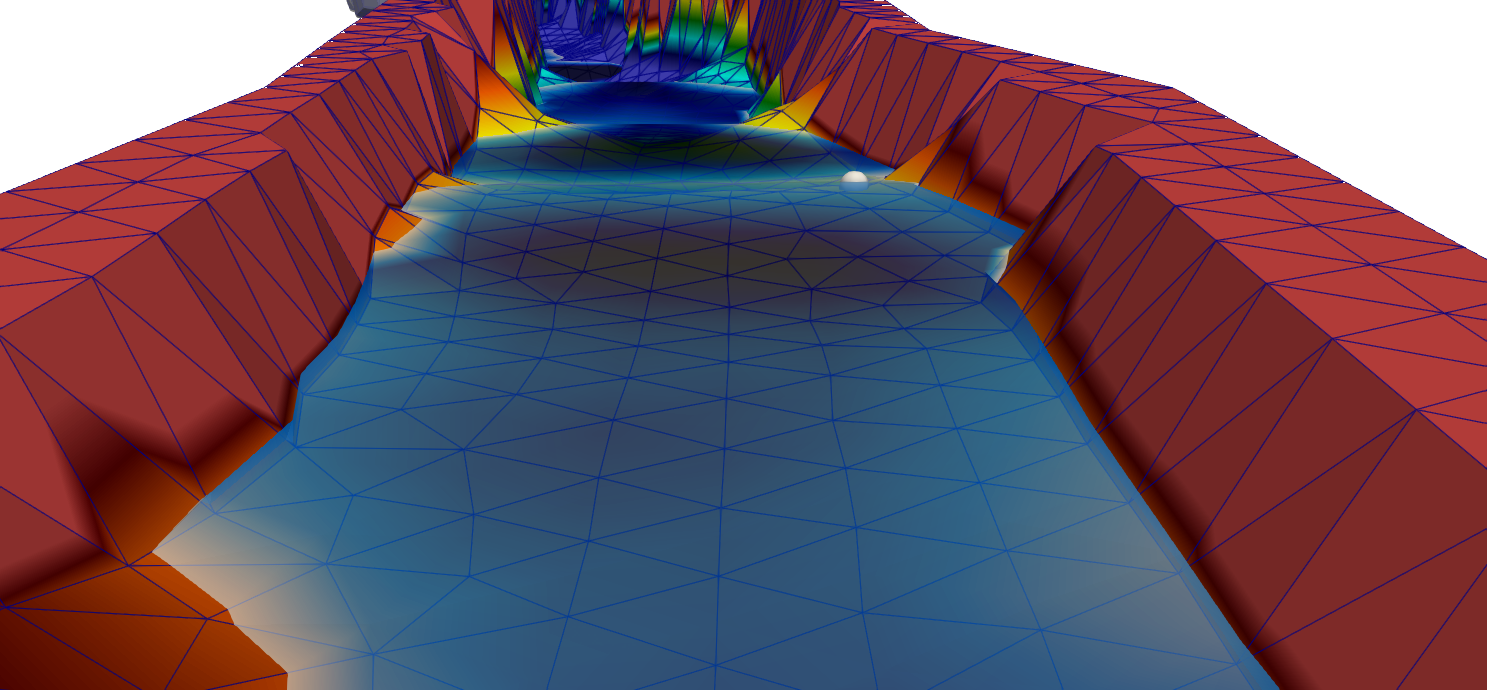}
         \caption{$t=4000$}
         \label{3d_result_t4000}
     \end{subfigure}
     \begin{subfigure}[b]{0.45\textwidth}
         \centering
         \includegraphics[width=\textwidth]{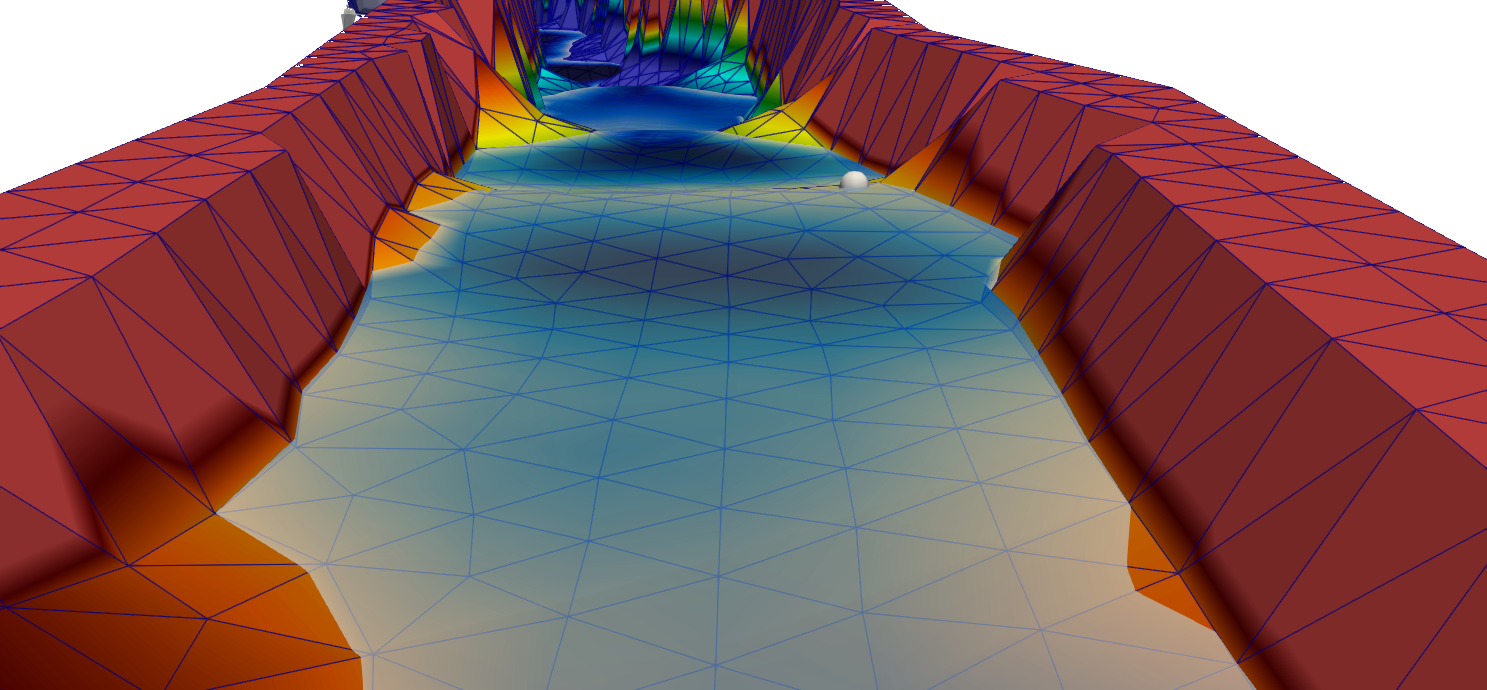}
         \caption{$t = 6000$}
         \label{3d_result_t6000}
     \end{subfigure}
     \hfill
     \begin{subfigure}[b]{0.45\textwidth}
         \centering
         \includegraphics[width=\textwidth]{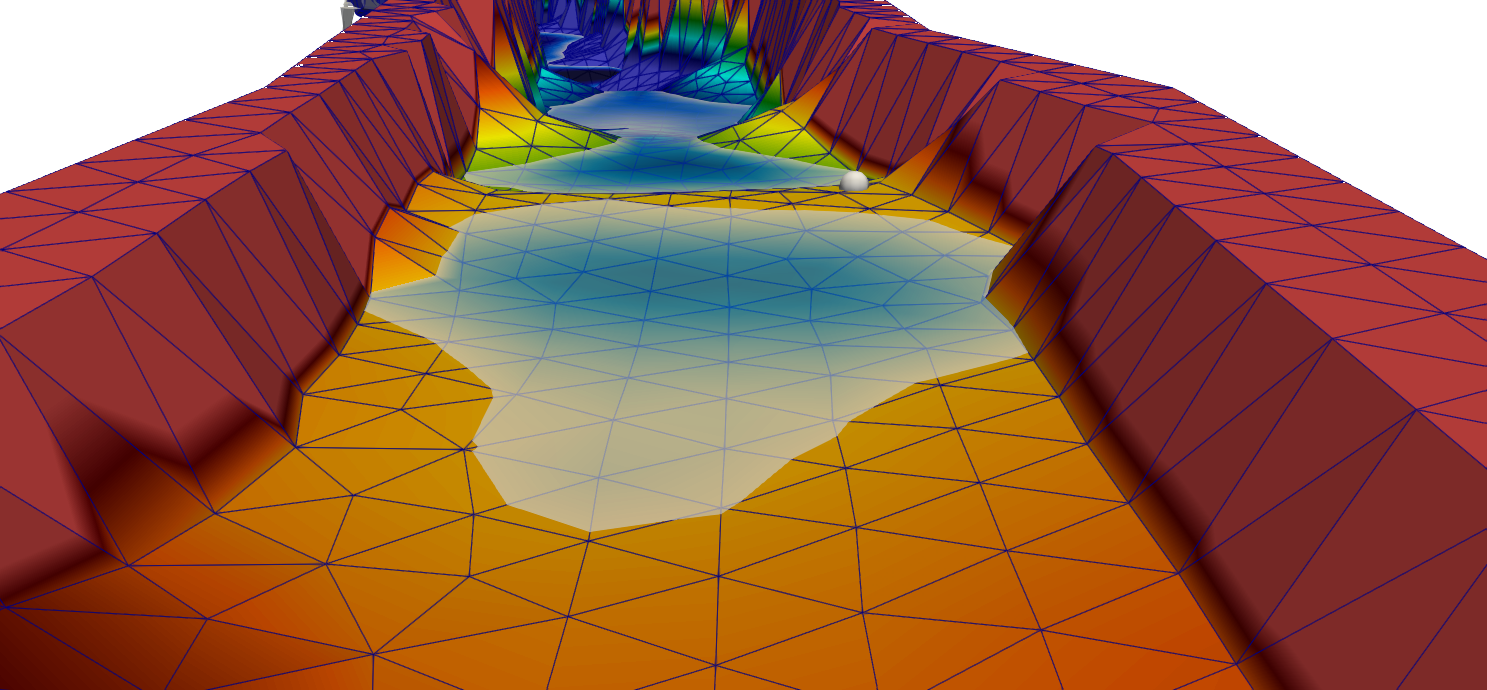}
         \caption{$t=8000$}
         \label{3d_result_t8000}
     \end{subfigure}
    \caption{3D visualization of the water flow looking down the river valley. The little gray sphere shows the location of measurement probe 6.}
    \label{probe6_3d}
\end{figure}

The hydrographs recorded at the seven probe locations are plotted in Fig.~\ref{fig:water_depth_test5}. Based on their behavior, the probes can be divided into three clear categories. The first category includes probes that retain residual water after the flood peak recedes (Probes 1, 3, and 4). These are all located in deeper sections of the main channel. The second category consists of probes that become completely dry by the end of the simulation (Probes 2, 6, and 7). This group includes probes on shallower ground (6 and 7) as well as a location within the main channel (2). The third category is unique to Probe 5, which is situated at the downstream domain outlet. It acts as a collection point, showing a continuous accumulation of water without a discernible peak. Note that probe 5 is also located in the shallow region. Note that subfigures in Fig.~\ref{fig:water_depth_test5} have different ranges on the vertical axis. 

The hydrographs show good agreement with the benchmark data from \cite{neelz2010benchmarking}, indicating that our proposed scheme, especially the cut-cell DG method, can cope with a realistic bathymmetry. Convergence is obtained at all probes with very good results for the CUDGM. Except for probe 6, already the coarsest level DG solution is very accurate.

The three-dimensional visualization of the water flow close to probe 6 demonstrates the complexity of the bathymmetry. The bathymmetry of that region is shown in Fig.~\ref{3d_result_t0}. The flow situation at probe 6 is difficult because it is located on an elevated hump and is also close to the valley boundary (the high slope region). The water initially enters the region at approximately $t = 1500 \text{ s}$ and subsequently passes through the probe location around $t = 2500 s$. Following this inflow, the water is gradually drained from the area. A small volume of water is observed flowing across the hump and still passing Probe 6, as shown in Fig.~\ref{3d_result_t4000}. Then, this flow detaches at approximately $t = 6000 s$ (Fig.~\ref{3d_result_t6000}), eventually separating entirely into two disconnected water bodies (Fig.~\ref{3d_result_t8000}).

\section{Conclusion} \label{conclusion section}

In this paper we presented a cut-cell upwind discontinuous Galerkin method (CUDGM) as well as a finite volume method on Voronoi cells (VFVM) for the diffusive wave approximation of the shallow water equations. Both methods maintain non-negative water heights and can handle unstructured triangular meshes. While the VFVM is restricted to Delauney triangulations and piecewise constant bathymmetry representation, the CUDGM in contrast can handle general unstructured meshes as well as continuous and discontinuous bathymmetries. For the Barenblatt solution on an inclined plane we demonstrate full second-order accuracy of the CUDGM including the free boundary. The VFVM on the other hand is only first-order accurate. Across three different numerical examples it is shown that the CUDGM achieves the same accuracy as the VFVM on three times coarser meshes. As a result, the method is also computationally more efficient.  

\section*{Acknowledgements}
The first author would like to acknowledge the financial support provided by the Development and Promotion of Science and Technology (DPST) Project, Thailand.
\bibliographystyle{plain} 
\bibliography{Reference}

\begin{thebibliography}{10}

\bibitem{alexander1977diagonally}
Roger Alexander.
\newblock Diagonally implicit runge--kutta methods for stiff ode’s.
\newblock {\em SIAM Journal on Numerical Analysis}, 14(6):1006--1021, 1977.

\bibitem{ALONSOSANTILLANADAWSON2008}
R.~Alonso, M.~Santillana, and C.~Dawson.
\newblock On the diffusive wave approximation of the shallow water equations.
\newblock {\em European Journal of Applied Mathematics}, 19(5):575–606, 2008.

\bibitem{Arico2018}
Costanza Aricò and Carmelo Nasello.
\newblock Comparative analyses between the zero-inertia and fully dynamic
  models of the shallow water equations for unsteady overland flow propagation.
\newblock {\em Water}, 10(1), 2018.

\bibitem{AUDUSSE2005311}
Emmanuel Audusse and Marie-Odile Bristeau.
\newblock A well-balanced positivity preserving “second-order” scheme for
  shallow water flows on unstructured meshes.
\newblock {\em Journal of Computational Physics}, 206(1):311--333, 2005.

\bibitem{BADIA2018533}
Santiago Badia, Francesc Verdugo, and Alberto~F. Martín.
\newblock The aggregated unfitted finite element method for elliptic problems.
\newblock {\em Computer Methods in Applied Mechanics and Engineering},
  336:533--553, 2018.

\bibitem{bastian2014fully}
Peter Bastian.
\newblock A fully-coupled discontinuous galerkin method for two-phase flow in
  porous media with discontinuous capillary pressure.
\newblock {\em Computational Geosciences}, 18(5):779--796, 2014.

\bibitem{Bastian2009}
Peter Bastian and Christian Engwer.
\newblock An unfitted finite element method using discontinuous galerkin.
\newblock {\em International Journal for Numerical Methods in Engineering},
  79(12):1557--1576, 2009.

\bibitem{bauer2006regional}
Peter Bauer, Thomas Gumbricht, and Wolfgang Kinzelbach.
\newblock A regional coupled surface water/groundwater model of the okavango
  delta, botswana.
\newblock {\em Water Resources Research}, 42(4), 2006.

\bibitem{bogelein2021existence}
Verena B{\"o}gelein, Nicolas Dietrich, and Matias Vestberg.
\newblock Existence of solutions to a diffusive shallow medium equation.
\newblock {\em Journal of Evolution Equations}, 21:845--889, 2021.

\bibitem{BOLSTER2002221}
Carl~H Bolster and James~E Saiers.
\newblock Development and evaluation of a mathematical model for surface-water
  flow within the shark river slough of the florida everglades.
\newblock {\em Journal of Hydrology}, 259(1):221--235, 2002.

\bibitem{BURMAN20101217}
Erik Burman.
\newblock Ghost penalty.
\newblock {\em Comptes Rendus Mathematique}, 348(21):1217--1220, 2010.

\bibitem{Burman2015}
Erik Burman, Susanne Claus, Peter Hansbo, Mats~G. Larson, and André Massing.
\newblock Cutfem: Discretizing geometry and partial differential equations.
\newblock {\em International Journal for Numerical Methods in Engineering},
  104(7):472--501, 2015.

\bibitem{caviedes2018cellular}
Daniel Caviedes-Voulli{\`e}me, Javier Fern{\'a}ndez-Pato, and Christoph Hinz.
\newblock Cellular automata and finite volume solvers converge for 2d shallow
  flow modelling for hydrological modelling.
\newblock {\em Journal of hydrology}, 563:411--417, 2018.

\bibitem{CAVIEDESVOULLIEME2020124663}
Daniel Caviedes-Voullième, Javier Fernández-Pato, and Christoph Hinz.
\newblock Performance assessment of 2d zero-inertia and shallow water models
  for simulating rainfall-runoff processes.
\newblock {\em Journal of Hydrology}, 584:124663, 2020.

\bibitem{cea2010experimental}
L~Cea, M~Garrido, and J~Puertas.
\newblock Experimental validation of two-dimensional depth-averaged models for
  forecasting rainfall--runoff from precipitation data in urban areas.
\newblock {\em Journal of Hydrology}, 382(1-4):88--102, 2010.

\bibitem{CEA201088}
L.~Cea, M.~Garrido, and J.~Puertas.
\newblock Experimental validation of two-dimensional depth-averaged models for
  forecasting rainfall–runoff from precipitation data in urban areas.
\newblock {\em Journal of Hydrology}, 382(1):88--102, 2010.

\bibitem{COSTABILE2017141}
Pierfranco Costabile, Carmelina Costanzo, and Francesco Macchione.
\newblock Performances and limitations of the diffusive approximation of the
  2-d shallow water equations for flood simulation in urban and rural areas.
\newblock {\em Applied Numerical Mathematics}, 116:141--156, 2017.
\newblock New Trends in Numerical Analysis: Theory, Methods, Algorithms and
  Applications (NETNA 2015).

\bibitem{costabile2017performances}
Pierfranco Costabile, Carmelina Costanzo, and Francesco Macchione.
\newblock Performances and limitations of the diffusive approximation of the
  2-d shallow water equations for flood simulation in urban and rural areas.
\newblock {\em Applied Numerical Mathematics}, 116:141--156, 2017.

\bibitem{de2019discontinuous}
Thomas De~Maet.
\newblock {\em Discontinuous Galerkin models for coupled surface-subsurface
  flow}.
\newblock PhD thesis, UCL-Universit{\'e} Catholique de Louvain, 2019.

\bibitem{DIETRICH201145}
J.C. Dietrich, M.~Zijlema, J.J. Westerink, L.H. Holthuijsen, C.~Dawson, R.A.
  Luettich, R.E. Jensen, J.M. Smith, G.S. Stelling, and G.W. Stone.
\newblock Modeling hurricane waves and storm surge using integrally-coupled,
  scalable computations.
\newblock {\em Coastal Engineering}, 58(1):45--65, 2011.

\bibitem{engwersisc2020}
Christian Engwer, Sandra May, Andreas N\"{u}\ss{}ing, and Florian
  Streitb\"{u}rger.
\newblock A stabilized dg cut cell method for discretizing the linear transport
  equation.
\newblock {\em SIAM Journal on Scientific Computing}, 42(6):A3677--A3703, 2020.

\bibitem{ern2009discontinuous}
Alexandre Ern, Annette~F Stephansen, and Paolo Zunino.
\newblock A discontinuous galerkin method with weighted averages for
  advection--diffusion equations with locally small and anisotropic
  diffusivity.
\newblock {\em IMA Journal of Numerical Analysis}, 29(2):235--256, 2009.

\bibitem{fernandez20162d}
J~Fern{\'a}ndez-Pato and P~Garc{\'\i}a-Navarro.
\newblock 2d zero-inertia model for solution of overland flow problems in
  flexible meshes.
\newblock {\em Journal of Hydrologic Engineering}, 21(11):04016038, 2016.

\bibitem{HANSBO20025537}
Anita Hansbo and Peter Hansbo.
\newblock An unfitted finite element method, based on {N}itsche’s method, for
  elliptic interface problems.
\newblock {\em Computer Methods in Applied Mechanics and Engineering},
  191(47):5537--5552, 2002.

\bibitem{Heimann2013}
F.~Heimann, C.~Engwer, O.~Ippisch, and P.~Bastian.
\newblock An unfitted interior penalty discontinuous galerkin method for
  incompressible navier–stokes two-phase flow.
\newblock {\em International Journal for Numerical Methods in Fluids},
  71(3):269--293, 2013.

\bibitem{HUNTER2007208}
Neil~M. Hunter, Paul~D. Bates, Matthew~S. Horritt, and Matthew~D. Wilson.
\newblock Simple spatially-distributed models for predicting flood inundation:
  A review.
\newblock {\em Geomorphology}, 90(3):208--225, 2007.
\newblock Reduced-Complexity Geomorphological Modelling for River and Catchment
  Management.

\bibitem{hunter2005adaptive}
Neil~M Hunter, Matthew~S Horritt, Paul~D Bates, Matthew~D Wilson, and Micha~GF
  Werner.
\newblock An adaptive time step solution for raster-based storage cell
  modelling of floodplain inundation.
\newblock {\em Advances in water resources}, 28(9):975--991, 2005.

\bibitem{liu2021parallel}
Shuang Liu, Qi~Tang, and Xian-zhu Tang.
\newblock A parallel cut-cell algorithm for the free-boundary grad--shafranov
  problem.
\newblock {\em SIAM Journal on Scientific Computing}, 43(6):B1198--B1225, 2021.

\bibitem{Martins2017}
Ricardo Martins, Jorge Leandro, Albert~S. Chen, and Slobodan Djordjević.
\newblock A comparison of three dual drainage models: shallow water vs local
  inertial vs diffusive wave.
\newblock {\em Journal of Hydroinformatics}, 19(3):331--348, 2017.

\bibitem{MIGNOT2006186}
E.~Mignot, A.~Paquier, and S.~Haider.
\newblock Modeling floods in a dense urban area using 2d shallow water
  equations.
\newblock {\em Journal of Hydrology}, 327(1):186--199, 2006.

\bibitem{neelz2010benchmarking}
Sylvain N{\'e}elz, Garry Pender, et~al.
\newblock Benchmarking of 2d hydraulic modelling packages.
\newblock Bristol: Environment Agency, 2010.

\bibitem{PaniconiPutti2015}
Claudio Paniconi and Mario Putti.
\newblock Physically based modeling in catchment hydrology at 50: Survey and
  outlook.
\newblock {\em Water Resources Research}, 51(9):7090--7129, 2015.

\bibitem{ponce1977shallow}
Victor~Miguel Ponce and Daryl~B Simons.
\newblock Shallow wave propagation in open channel flow.
\newblock {\em Journal of the Hydraulics Division}, 103(12):1461--1476, 1977.

\bibitem{riviere2008discontinuous}
B{\'e}atrice Rivi{\`e}re.
\newblock {\em Discontinuous Galerkin methods for solving elliptic and
  parabolic equations: theory and implementation}.
\newblock SIAM, 2008.

\bibitem{santillana2010local}
Mauricio Santillana and Clint Dawson.
\newblock A local discontinuous galerkin method for a doubly nonlinear
  diffusion equation arising in shallow water modeling.
\newblock {\em Computer methods in applied mechanics and engineering},
  199(23-24):1424--1436, 2010.

\bibitem{santillana2010numerical}
Mauricio Santillana and Clint Dawson.
\newblock A numerical approach to study the properties of solutions of the
  diffusive wave approximation of the shallow water equations.
\newblock {\em Computational Geosciences}, 14:31--53, 2010.

\bibitem{savant2019urban}
Gaurav Savant, Corey~J. Trahan, Lucas Pettey, Tate~O. McAlpin, Gary~L. Bell,
  and Charles~J. McKnight.
\newblock Urban and overland flow modeling with dynamic adaptive mesh and
  implicit diffusive wave equation solver.
\newblock {\em Journal of Hydrology}, 573:13--30, 2019.

\bibitem{schneiders2013accurate}
Lennart Schneiders, Daniel Hartmann, Matthias Meinke, and Wolfgang
  Schr{\"o}der.
\newblock An accurate moving boundary formulation in cut-cell methods.
\newblock {\em Journal of Computational Physics}, 235:786--809, 2013.

\bibitem{singer2019local}
Thomas Singer and Matias Vestberg.
\newblock Local boundedness of weak solutions to the diffusive wave
  approximation of the shallow water equations.
\newblock {\em Journal of Differential Equations}, 266(6):3014--3033, 2019.

\bibitem{SOHEB2024131063}
Mohd Soheb, Peter Bastian, Susanne Schmidt, Shaktiman Singh, Himanshu Kaushik,
  Alagappan Ramanathan, and Marcus Nüsser.
\newblock Surface and subsurface flow of a glacierised catchment in the
  cold-arid region of ladakh, trans-himalaya.
\newblock {\em Journal of Hydrology}, 635:131063, 2024.

\bibitem{SZYMKIEWICZ2012165}
Romuald Szymkiewicz and Dariusz Gasiorowski.
\newblock Simulation of unsteady flow over floodplain using the diffusive wave
  equation and the modified finite element method.
\newblock {\em Journal of Hydrology}, 464-465:165--175, 2012.

\bibitem{szymkiewicz2012simulation}
Romuald Szymkiewicz and Dariusz Gasiorowski.
\newblock Simulation of unsteady flow over floodplain using the diffusive wave
  equation and the modified finite element method.
\newblock {\em Journal of hydrology}, 464:165--175, 2012.

\bibitem{tan1992shallow}
Wei-Yan Tan.
\newblock {\em Shallow water hydrodynamics: Mathematical theory and numerical
  solution for a two-dimensional system of shallow-water equations}, volume~55.
\newblock Elsevier, 1992.

\bibitem{Thornton2022}
J.~M. Thornton, R.~Therrien, G.~Mariéthoz, N.~Linde, and P.~Brunner.
\newblock Simulating fully-integrated hydrological dynamics in complex alpine
  headwaters: Potential and challenges.
\newblock {\em Water Resources Research}, 58(4):e2020WR029390, 2022.

\bibitem{vantiel2024}
Marit van Tiel, Caroline Aubry-Wake, Lauren Somers, Christoff Andermann,
  Francesco Avanzi, Michel Baraer, Gabriele Chiogna, Cl{\'e}mence Daigre,
  Soumik Das, Fabian Drenkhan, Daniel Farinotti, Catriona~L. Fyffe, Inge
  de~Graaf, Sarah Hanus, Walter Immerzeel, Franziska Koch, Jeffrey~M. McKenzie,
  Tom M{\"u}ller, Andrea~L. Popp, Zarina Saidaliyeva, Bettina Schaefli,
  Oliver~S. Schilling, Kapiolani Teagai, James~M. Thornton, and Vadim Yapiyev.
\newblock Cryosphere--groundwater connectivity is a missing link in the
  mountain water cycle.
\newblock {\em Nature Water}, 2(7):624--637, 2024.

\bibitem{VINCENT201959}
Aude Vincent, Sophie Violette, and Gudfinna Adalgeirsdóttir.
\newblock Groundwater in catchments headed by temperate glaciers: A review.
\newblock {\em Earth-Science Reviews}, 188:59--76, 2019.

\bibitem{wang2011positivity}
Yueling Wang, Qiuhua Liang, Georges Kesserwani, and Jim~W Hall.
\newblock A positivity-preserving zero-inertia model for flood simulation.
\newblock {\em Computers \& Fluids}, 46(1):505--511, 2011.

\bibitem{Wood2011}
Eric~F. Wood, Joshua~K. Roundy, Tara~J. Troy, L.~P.~H. van Beek, Marc F.~P.
  Bierkens, Eleanor Blyth, Ad~de~Roo, Petra Döll, Mike Ek, James Famiglietti,
  David Gochis, Nick van~de Giesen, Paul Houser, Peter~R. Jaffé, Stefan
  Kollet, Bernhard Lehner, Dennis~P. Lettenmaier, Christa Peters-Lidard,
  Murugesu Sivapalan, Justin Sheffield, Andrew Wade, and Paul Whitehead.
\newblock Hyperresolution global land surface modeling: Meeting a grand
  challenge for monitoring earth's terrestrial water.
\newblock {\em Water Resources Research}, 47(5), 2011.

\end{thebibliography}

\end{document}